\documentclass{siamltex}
\usepackage{graphicx}
\graphicspath{{/home/bette/Dropbox/Research/RBF-PU-LS/}
              {/home/bette/Dropbox/Research/RBF-PU-LS/NewFigs/}}
\usepackage{amsmath}
\usepackage{amssymb}

\usepackage{color}

\newcommand{\ep}[1]{\ensuremath{\varepsilon^{#1}}}
\newcommand{\x}{\ensuremath{\underline{x}}}
\newcommand{\y}{\ensuremath{\underline{y}}}

\newcommand{\ut}{\ensuremath{\tilde{u}}}

\newtheorem{assume}[theorem]{Assumption}
\newtheorem{estimate}[theorem]{Estimate}

\begin{document}
\title{A least squares radial basis function partition of unity method for solving PDEs}
\author{Elisabeth Larsson\thanks{Uppsala University, Department of Information Technology, Box~337, SE-751~05 Uppsala, Sweden
(elisabeth.larsson@it.uu.se, victor.shcherbakov@it.uu.se).}
\and Victor Shcherbakov\footnotemark[1]
\and Alfa Heryudono\thanks{University of Massachusetts Dartmouth, Department of Mathematics, 285 Old Westport Road, Dartmouth, MA 02747, USA (aheryudono@umassd.edu). The work was supported by the European Commission CORDIS Marie Curie FP7 program Grant \#235730 and National Science Foundation DMS Grant \#1552238.}
}
\maketitle

\begin{abstract}
Recently, collocation based radial basis function (RBF) partition of unity methods (PUM) for solving partial differential equations have been formulated and investigated numerically and theoretically. When combined with stable evaluation methods such as the RBF-QR method, high order convergence rates can be achieved and sustained under refinement. However, some numerical issues remain. The method is sensitive to the node layout, and condition numbers increase with the refinement level. Here, we propose a modified formulation based on least squares approximation. We show that the sensitivity to node layout is removed and that conditioning can be controlled through oversampling. We derive theoretical error estimates both for the collocation and least squares RBF-PUM. Numerical experiments are performed for the Poisson equation in two and three space dimensions for regular and irregular geometries. The convergence experiments confirm the theoretical estimates, and the least squares formulation is shown to be 5--10 times faster than the collocation formulation for the same accuracy. 

\end{abstract}

\pagestyle{myheadings}
\thispagestyle{plain}
\markboth{E. LARSSON, V. SHCHERBAKOV, A. HERYUDONO}{A LEAST SQUARES RBF-PUM}

\section{Introduction}
Radial basis function (RBF) approximation methods for partial differential equations (PDEs) have several important advantages. For PDEs with smooth solutions, approximation with infinitely smooth RBFs provides spectral convergence properties for non-trivial geometries~\cite{Schaback07,RieZwi10,RieZwi14}. It is easy to formulate and implement RBF methods in any number of dimensions due to the meshfree nature of the methods and the reduction of geometrical properties to computations of pairwise distances. A global RBF approximation $\ut(\x)$ to a function $u(\x)$ has the form
\begin{equation}
\tilde{u}(\x)=\sum_{i=1}^{n}\lambda_{i}\phi(\|\x-\x_{i}\|),
\label{eq:rbf1}
\end{equation}
where $\x\in\mathbb{R}^d$, $\phi(r)$ is an RBF, $\x_1,\ldots,\x_n$ are scattered node points at which the RBFs are centered, and $\lambda_i$ are coefficients to be determined from given data. Whether an interpolation problem or a PDE problem is solved using collocation~\cite{LaFo03,FlyWri09}, least squares approximation~\cite{KwoLi09,AlWaJeRi12}, or a Galerkin approach~\cite{Wend99,KoLa13}, the resulting linear systems that need to be solved are dense. Due to the high order convergence rate, the systems are comparatively small in size, but for geometrically large scale problems in more than two space dimensions, both the computational cost and the storage requirements still become prohibitive.

Two main directions of research with the purpose of reducing the computational cost of RBF methods through localization are currently pursued. The first one, RBF-FD methods, can be seen as a generalization of finite difference methods, but with stencils supported on scattered node sets. To our knowledge, Tolstykh~\cite{Tol00} was first to publish the method, which since has been extensively researched, see, e.g.,~\cite{ShuDiYe03,WriFo06,FoLeh11,DavOanh11,FBWStC12,DavScha16}. The current state of the art is well presented in~\cite{FoFly15,FFBB16,BaFlyFoBa17}.

The other direction, which is the focus of this paper is RBF-based partition of unity methods (RBF-PUM). The idea of combining RBF approximations with partition of unity was suggested already by Babu\v{s}ka and Melenk~\cite{BaMe97}. RBF-PUM was explored for interpolation purposes in combination with compactly supported RBFs by Wendland~\cite{Wend02}, and further discussed in the book on meshfree approximation by Fasshauer~\cite{Fass07}. Lately, Cavoretto, De Rossi et al.~have explored various method and implementation aspects of RBF-PUM for interpolation of non-uniform scattered data~\cite{CaDeRo12,CaDeRo14,CaDeRo15,CaDeRoPe16,CaDeRoPe16b}.

In the forthcoming paper~\cite{LaHer17} (see also~\cite{LaHer12,SVHL15}), we derive a collocation based RBF-PUM (C-RBF-PUM) for PDE problems, and provide theoretical results for the approximation errors of such an approach. The approximation error drives the convergence of a PDE solution~\cite{Schaback07}, but does not take issues related to well-posedness and conditioning into account. The collocation method works well, and has been used successfully for option pricing problems (parabolic PDEs)~\cite{SVHL15,ShchLa16,Shch16,BENCHOP} as well as for glacier modeling~\cite{AhlShch16}. A key to the success of the method is the use of the RBF-QR method for stable evaluation~\cite{FoPi07,FoLaFly11,LLHF13}. Despite the overall positive results, there are some issues to consider: The method exhibits some sensitivity to the node layout, especially near boundaries where it is difficult to maintain a quasi uniform node structure, and the linear systems become increasingly ill-conditioned when the problem size grows, making it practically difficult to address large scale problems.

In this paper, we move away from the collocation approach in favor of a least squares approach. We allow the node points to be decoupled from the problem geometry, thus simplifying node generation while allowing for high quality node layouts. The features of the geometry are instead captured by the choice of the least squares evaluation points that are used to enforce the PDE and its boundary conditions. The oversampling, resulting in an overdetermined linear system that is solved using least squares, removes the robustness issues related both to boundaries and problem scale. We derive full error estimates for elliptic PDEs for both C-RBF-PUM and the least squares (LS-RBF-PUM) approach. Furthermore, we perform extensive numerical experiments for elliptic PDEs in two and three space dimensions to illustrate the significantly improved properties of the new formulation of the method.

The paper is organized as follows: In Section~\ref{sec:model}, the Poisson test problems are discussed. Section~\ref{sec:rbfpum} derives the two RBF-PUM approaches, and in Section~\ref{sec:theory} theoretical convergence estimates are provided. Numerical experiments on convergence, robustness, and computational cost are shown in Section~\ref{sec:exp} for two-dimensional and three-dimensional problems. The final section in the paper contains a discussion of the methods and results.




\section{The Poisson test problems}\label{sec:model}

We have chosen to use the linear, elliptic, time-independent Poisson equation, with Dirichlet boundary conditions, in two and three spatial dimensions as test problems to compare the two RBF-PUM formulations that are investigated in the paper. With this choice, we focus solely on the spatial PDE approximation properties and avoid complications arising from an additional time discretization. Achieving competitive performance for the Poisson equation is a requirement for later moving to more advanced PDEs. The problem in its general form is 
\begin{equation}
\left\{\begin{array}{rcll}
-\Delta u(\x) & = & f(\x), & \x\in\Omega,\\
u(\x) & = & g(\x), & \x\in\partial\Omega,
\end{array}\right.
\label{eq:poisson}
\end{equation}
where $\x=(x_1,\ldots,x_d)\in\mathbb{R}^d$. 
We need some general assumptions to hold for the geometry of the domain $\Omega$ to be able to later derive convergence estimates for RBF approximations.
\begin{assume}
\label{ass:cone}
The domain $\Omega \subset \mathbb{R}^d$ is an open, bounded domain with Lipschitz boundary, that satisfies an interior cone condition~\cite{Wend05} with maximum radius $\mathcal{R}$ and angle $\nu$.
\end{assume}
%

For the numerical experiments in $\mathbb{R}^2$, we use three different domains. The box
\begin{equation}
\Omega_B=\{\x : |x_i|\leq 2,\ i=1,2\},
\end{equation}
a star-shaped, non-convex domain with smooth boundary, defined using polar coordinates as
\begin{equation}
\Omega_S =\{\x=(r,\theta) : r\leq  2(0.7+0.12(\sin(6\theta)+\sin(3\theta))),\ \theta\in[0,2\pi)\}, 
\end{equation}
and a polygonal non-convex domain with a Lipschitz boundary $\partial\Omega_L$ representing the mainland border of Sweden scaled to height 2,
\begin{equation}
\Omega_L =\{\x : \x \mbox{ inside } \partial\Omega_L\}.
\end{equation}
As an example, the interior cone condition holds with $\nu=\pi/4$ and $\mathcal{R}=1$ for $\Omega_B$.

For the experiments in $\mathbb{R}^3$, we have used two star-shaped domains, the unit sphere
\begin{equation}
\Omega_U=\{\x : \|\x\|\leq 1\},
\end{equation}
for which the interior cone condition holds with $\nu=\pi/3$ and $\mathcal{R}=1$, and the non-convex domain
\begin{equation}
\Omega_Q=\{\x=(r,\theta,\varphi) : r\leq r_Q(\theta,\varphi),\ \theta\in[0,2\pi),\ \varphi\in[0,\pi]\},
\end{equation}
where $r_Q=\left(1+\sin^2(2\sin\varphi\cos\theta)\sin^2(2\sin\varphi\sin\theta)\sin^2(2\cos\varphi)\right)^{1/2}$.

We also make assumptions on the types of solutions we consider for approximation by smooth RBFs. In order to achieve high order convergence, the regularity of $u$ needs to be higher than what is strictly required by the problem itself.
\begin{assume}
\label{ass:sol}
The solution $u(\x)$ to~\eqref{eq:poisson} as well as its first and second order derivatives are bounded and the following holds
\begin{equation*}
u(\x)\in W_\infty^m(\Omega)\subset W_\infty^2(\Omega),
\end{equation*}
where $W_p^k(\Omega)$ is a Sobolev space and $u\in W_p^k(\Omega)$ implies that
$\displaystyle\frac{\partial^{|\alpha|} u}{\partial \x^\alpha}\in L_p(\Omega)$, $\forall \alpha : |\alpha|\leq k$.
\end{assume}



\begin{figure}[!htb]
\centering
\begin{tabular}{@{}cc@{}}
\includegraphics[width=0.49\textwidth]{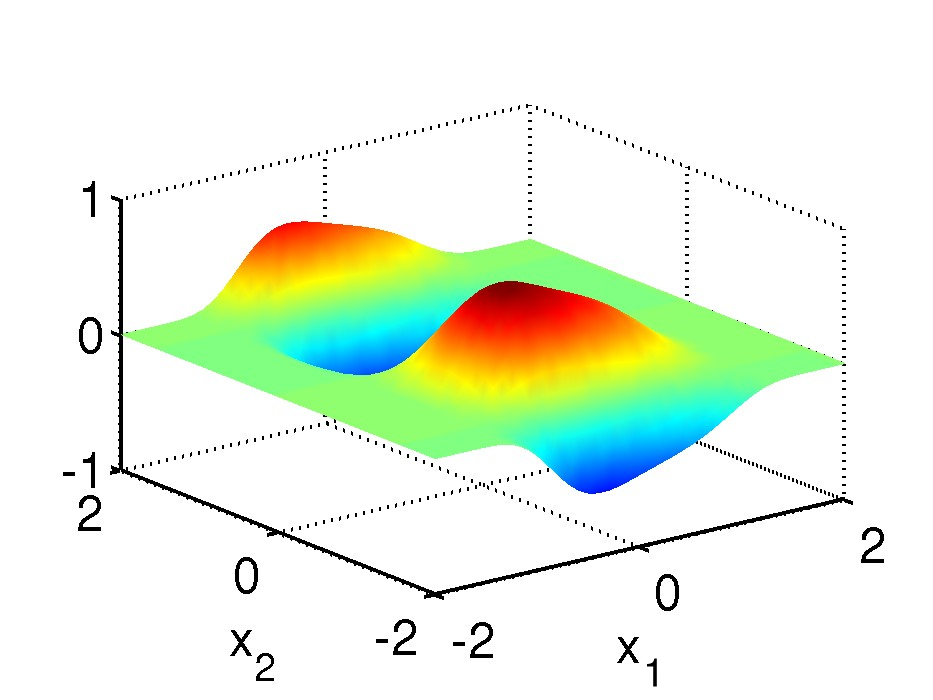} &
\includegraphics[width=0.45\textwidth]{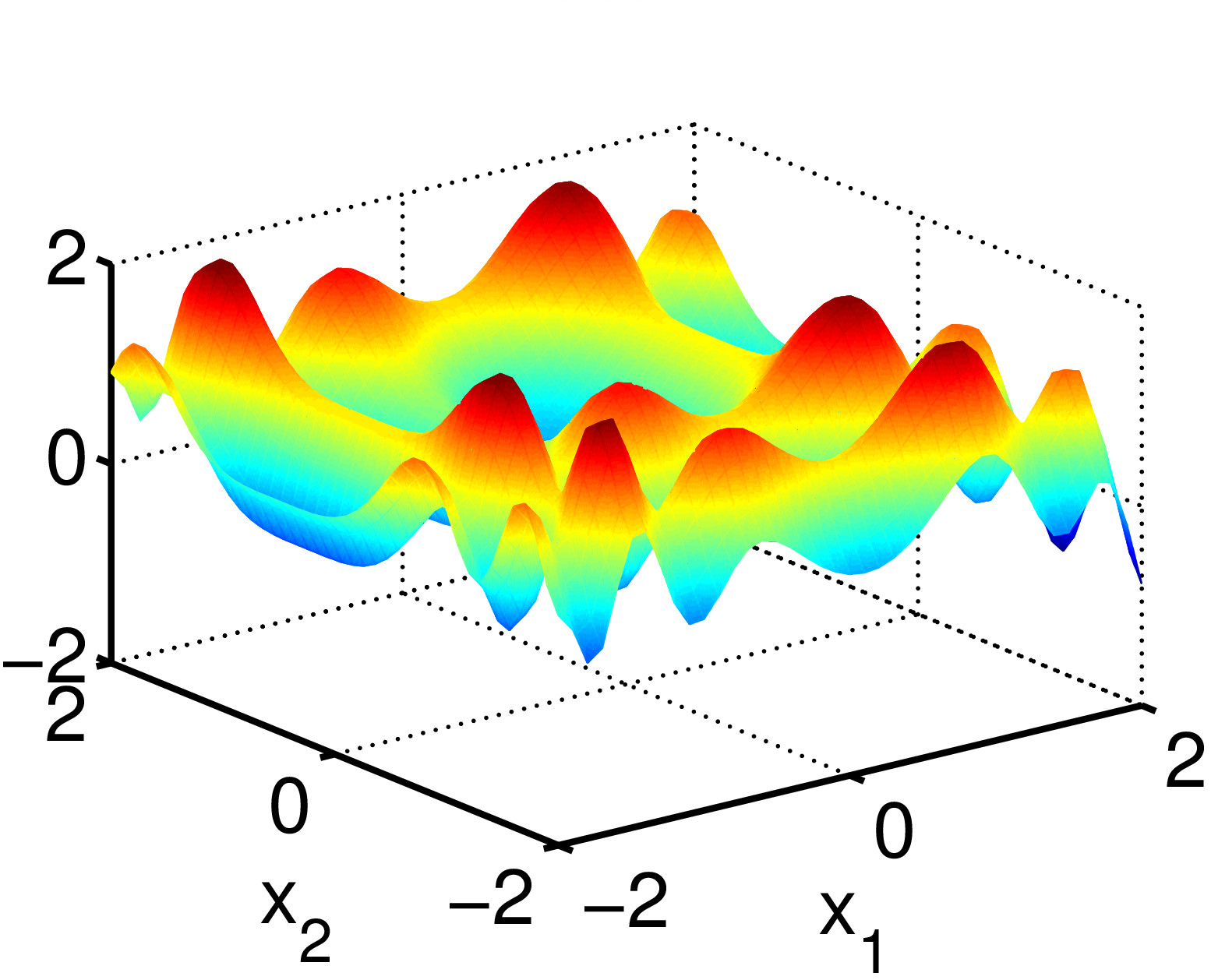}\\
\includegraphics[width=0.45\textwidth]{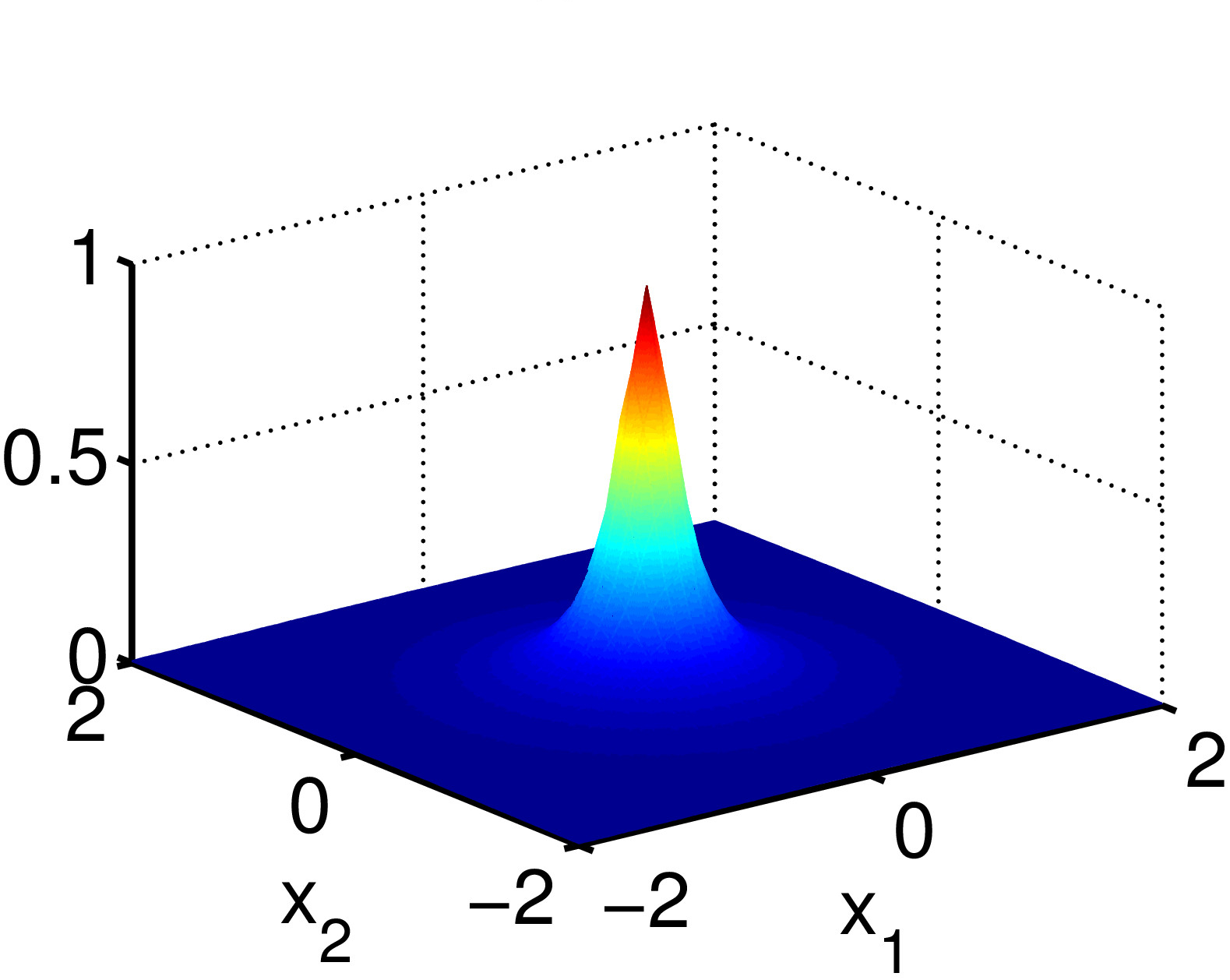} &
\includegraphics[width=0.45\textwidth]{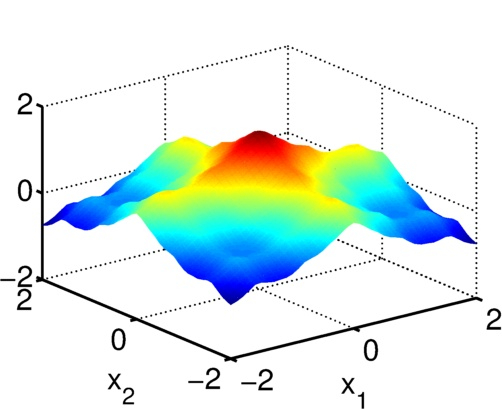}
\end{tabular}
\caption{The solution functions used in $\mathbb{R}^2$ displayed over the domain $\Omega_B$. The hyperbolic sine function $u_1$ (top left), the trigonometric combination $u_2$ (top right), the Runge function $u_3$ (bottom left), and the truncated non-analytic sum $u_4$ (bottom right).}
\label{fig:func}
\end{figure}
\begin{figure}[!htb]
\centering
\includegraphics[width=0.35\textwidth]{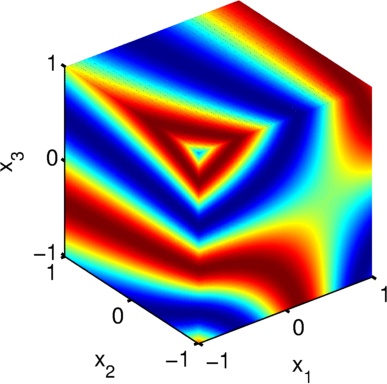}
\caption{The solution function $u_5$ used in $\mathbb{R}^3$ displayed over the unit cube. Function values range from -1 to 1.}
\label{fig:func3}
\end{figure}

For the numerical experiments, we use four different solution functions in $\mathbb{R}^2$ and one in $\mathbb{R}^3$ all in $W_\infty^\infty(\Omega)$ from which $f$ and $g$ in~\eqref{eq:poisson} are derived when solving the Poisson problem. The functions in $\mathbb{R}^2$ are illustrated in Figure~\ref{fig:func}, and are chosen to be increasingly difficult to approximate. The first function is a hyperbolic sine with a low number of oscillations within the domain,
%
\begin{equation}
u_1(\x)=\sinh\left(0.3(x_1-2)\sin(2x_2)\exp\left(-(x_1-0.1)^4\right)\right).
\end{equation}
The second function is a more oscillatory trigonometric combination
\begin{equation}
u_2(\x)=\sin\left(2(x_1-0.1)^2\right)\cos\left((x_1-0.3)^2\right)+\sin^2\left((x_2-0.5)^2\right).
\end{equation}
The third function is a sharp Runge type function. It is also equivalent to an inverse quadratic RBF with shape parameter $\ep{}=5$ placed at the origin,
\begin{equation}
u_3(\x)=\frac{1}{25x_1^2 + 25x_2^2 + 1}.
\end{equation}
The final $\mathbb{R}^2$ function is the first six modes of an expansion that for $j\rightarrow\infty$ becomes non-analytic,
\begin{equation}
u_4(\x)=\sum_{j=0}^{5}\exp(-\sqrt{2^j})(\cos(2^jx_1)+\cos(2^jx_2)).
%
\end{equation}
The $\mathbb{R}^3$ function we have used is given by
\begin{equation}
u_5(\x)=\sin\left(\frac{\pi(x_1-0.5)x_3}{\log(x_2+3)}\right),
\end{equation}
and is illustrated in Figure~\ref{fig:func3}.

For the theoretical convergence estimates derived in Section~\ref{sec:theory}, we need a well-posedness estimate that relates the norm of the solution $u$ of~\eqref{eq:poisson} to the data $f$ and $g$. For the case $f\equiv 0$, problem~\eqref{eq:poisson} is reduced to the Laplace equation, and the maximum principle holds for the solution $u_g$
\begin{equation}
\|u_g\|_{L_\infty(\Omega)}\leq \|g\|_{L_\infty(\partial\Omega)}.
\label{eq:max}
\end{equation}
If we instead have $g\equiv 0$, $\Omega$ satisfies Assumption~\ref{ass:cone}, and the solution $u_f\in W_p^1$ the classical Poincar\'e inequality holds
\begin{equation}
\|u_f\|_{L_p(\Omega)}\leq C\|\nabla u_f\|_{L_p(\Omega)}, \quad 1\leq p<\infty,
\label{eq:poincare}
\end{equation}
where the constant $C$ depends on $p$ and $\Omega$. Using the PDE~\eqref{eq:poisson} for $u_f\in W_2^2$, we can also through integration by parts and the Cauchy inequality derive
\begin{equation}
-\int_\Omega  u_f\Delta u_f=\int_\Omega \nabla u_f\cdot \nabla u_f=\|\nabla u_f\|_{L_2(\Omega)}^2=\int_\Omega  u_f f\leq \|u_f\|_{L_2(\Omega)}\|f\|_{L_2(\Omega)}.
\label{eq:grad}
\end{equation}
By combining~\eqref{eq:poincare} and~\eqref{eq:grad} (or applying the Poincar\'e inequality twice) we get
\begin{equation}
\|u_f\|_{L_2(\Omega)}\leq C^2\|f\|_{L_2(\Omega)}.
\label{eq:uf}
\end{equation}
By relating the $L_2$ and $L_\infty$ norms through
\[\|u\|_{L_2(\Omega)}^2=\int_\Omega|u|^2 \leq\max_\Omega|u|^2\int_\Omega1=|\Omega|\|u\|_{L_\infty(\Omega)}^2,\]
we can combine~\eqref{eq:max} and~\eqref{eq:uf} to arrive at the estimate
\begin{align}
\|u\|_{L_2(\Omega)} & =\|u_g+u_f\|_{L_2(\Omega)}\leq \|u_g\|_{L_2(\Omega)}+\|u_f\|_{L_2(\Omega)}\nonumber\\
&\leq \sqrt{|\Omega|}\|u_g\|_{L_\infty(\Omega)}+C^2\|f\|_{L_2(\Omega)}\nonumber\\
&\leq 
\sqrt{|\Omega|}\|g\|_{L_\infty(\partial\Omega)}+C^2\sqrt{|\Omega|}\|f\|_{L_\infty(\Omega)}.
\label{eq:est0}
\end{align}
Following~\cite{Schaback07}, to simplify notation, we define the operator $\mathcal{L}$ such that
\begin{equation}
\mathcal{L}u(\x) = \left\{\begin{array}{rl}
-\Delta u(\x),  & \x\in\Omega,\\
u(\x),  & \x\in\partial\Omega,
\end{array}\right.  
\end{equation}
and introduce a data norm defined as
\begin{equation}
\|u\|_F = \max(\|\Delta u\|_{L_\infty(\Omega)},\|u\|_{L_\infty(\partial\Omega)}).
\end{equation}
Then we can summarize the estimate~\eqref{eq:est0} and the corresponding assumptions as
\begin{estimate}
For a solution $u$ to the problem~\eqref{eq:poisson}, that satisfies Assumption~\ref{ass:sol}, over a domain $\Omega$ that satisfies Assumption~\ref{ass:cone}, it holds
\[\|u\|_{L_2(\Omega)}\leq C_P\|u\|_F,\]
where $C_P$ is a constant that depends on the shape and size of $\Omega$. 
\label{est:well}
\end{estimate}

\section{The RBF-based partition of unity methods}\label{sec:rbfpum}
In this section, we first provide a general description of partition of unity methods, then we discuss the local RBF approximations, and finally we combine these elements into the two different RBF-PUM formulations.
\subsection{Partition of unity methods}
\begin{figure}[!htb]
\centering
\includegraphics[width=0.35\textwidth]{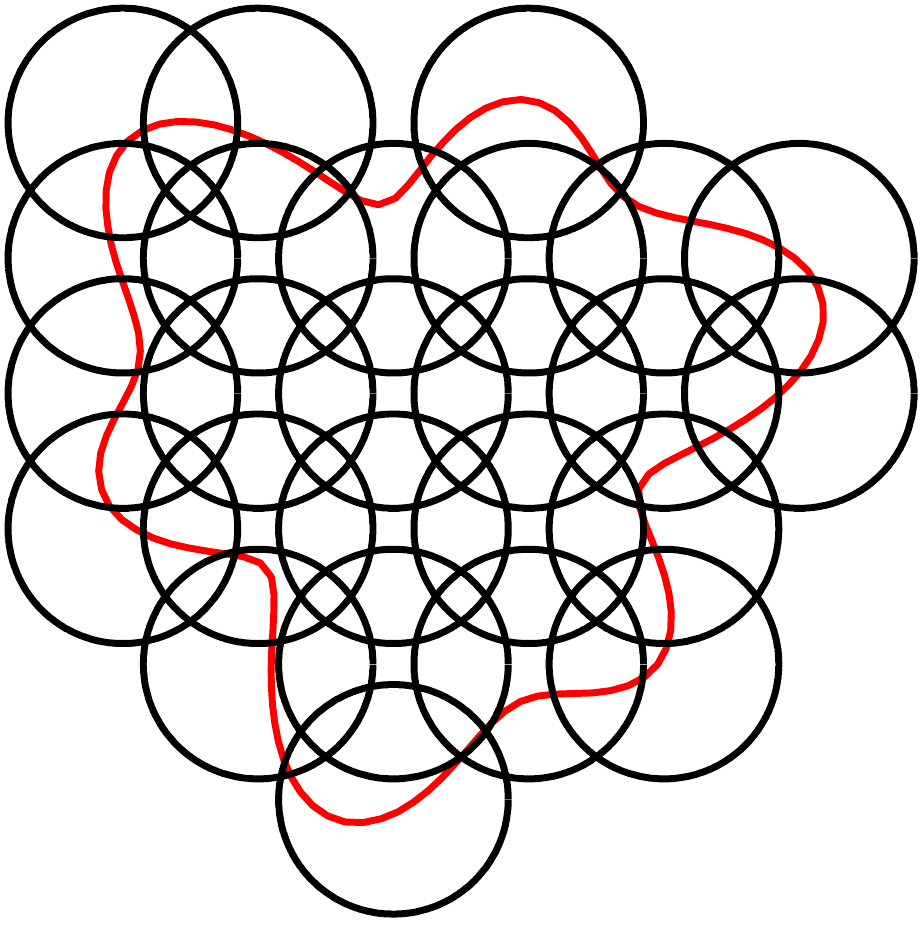}
\includegraphics[width=0.3\textwidth]{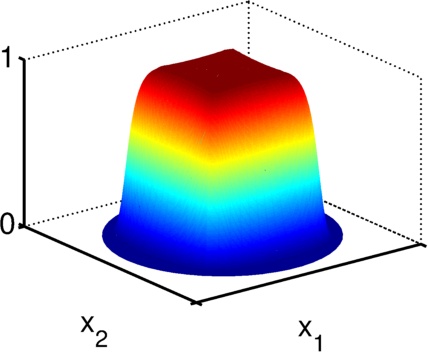}
\caption{To the left, the red curve is the outline of the domain $\Omega_S$ and the black circles are the boundaries of the overlapping circular patches $\Omega_j$, $j=1,\ldots,P$. To the right, a partition of unity weight function $w_j$ for one of the interior patches is shown.}
\label{fig:dom}
\end{figure}

To define a partition of unity method~\cite{BaMe97} for problem~\eqref{eq:poisson}, we construct a set of overlapping patches $\Omega_j$, $j=1,\ldots,P$ that form an open cover of the domain $\Omega$,
\[\bigcup_{j=1}^P\Omega_j\supseteq \Omega.\]
The amount of overlap between patches should be limited such that at most $K$ patches overlap at any given point. Throughout this paper, we choose to define the patches as discs in $\mathbb{R}^2$ and spheres in $\mathbb{R}^3$. Furthermore, the patch centers are chosen as vertices in an underlying Cartesian grid. With this choice, we can guarantee that the domain is covered and regulate the amount of overlap. This approach leads to $K=2^d$ in $\mathbb{R}^d$.
Other types of patches such as squares and cubes or ellipses~\cite{SVHL15} and ellipsoids can also be used, as well as less structured patch layouts~\cite{HLRvS16}.
An example of a cover of the star-shaped domain $\Omega_S$ with circular patches is shown in Figure~\ref{fig:dom}. Patches that do not contribute uniquely to the cover are pruned from the initial set that contains all patches $\Omega_j$ that have an intersection with $\Omega$, see~\cite{LaHer17} for a more detailed description.

In addition to the patches, we also construct partition of unity weight functions $w_j(\x)$, $j=1,\ldots,P$, subordinate to the open cover, such that
\[\sum_{j=1}^Pw_j(\x)=1,\quad \forall \x\in\Omega.\]
The weight function $w_j(\x)$ is compactly supported on $\Omega_j$, and
can be constructed by applying Shepard's method~\cite{Shepard68} to compactly supported generating functions $\varphi_j(\x)$, $j=1,\ldots,P$,
\begin{equation}
w_j(\x)=\frac{\varphi_j(\x)}{\sum_{i=1}^P\varphi_i(\x)},\qquad j=1,\ldots,P.
\end{equation}
The generating function needs to be smooth enough to support the differential operators of the problem to be solved. In our case, we have a second order elliptic PDE in strong form. We choose the generating function to be a compactly supported $C^2$ Wendland functions~\cite{Wendland95}, such as
\[\varphi(r)=(4r+1)(1-r^4)_+,\]
which is $C^2$ in up to three space dimensions.
In order to map the generating function to the patch $\Omega_j$ with center point $\underline{c}_j$ and radius $\rho_j$ we shift and scale the argument such that
\begin{equation}
\varphi_j(\x)=\varphi\left(\frac{\|\x-\underline{c}_j\|}{\rho_j}\right).
\end{equation}

A global partition of unity approximation $\ut(\x)$ to a function $u(\x)$ over $\Omega$ is formed as a weighted sum of local approximations $\ut_j(\x)$ on $\Omega_j$, using the partition of unity weight functions,
\begin{equation}
\ut(\x)=\sum_{j=1}^Pw_j(\x)\ut_j(\x).
\label{eq:PUM}
\end{equation} 
Partition of unity methods offer flexibility in the sense that the local approximations can be modified independently to match the local properties of the solution. In this article, we consider problems with smooth solutions, motivating the use of local smooth RBF approximations.

\subsection{RBF approximations and differentiation matrices}
We consider one of the local approximations $\ut_j(\x)$ on the patch $\Omega_j$.
RBF methods are meshfree, and approximations are defined on scattered node sets. We define two different scattered node sets, $X_j=\{\x_i^j\}_{i=1}^{n_j}$, at which the individual RBFs are centered, and $Y_j=\{\underline{y}_i^j\}_{i=1}^{m_j}$, where the RBF approximation is evaluated. Then we introduce a positive definite RBF $\phi(r)$, such as the Gaussian $\phi(r)=\exp(-\ep{2}r^2)$, or a conditionally positive definite RBF such as the multiquadric $\phi(r)=\sqrt{1+\ep{2}r^2}$, where $\ep{}$ is a shape parameter that determines the flatness of the basis functions. The standard form of the RBF approximation $\tilde{u}_j(\x)$ using basis functions centered at $X_j$ is 
\begin{equation}
\tilde{u}_j(\x)=\sum_{i=1}^{n_j}\lambda_{i}^j\phi(\|\x-\x_{i}^j\|)
\label{eq:RBF}
\end{equation}
where $\lambda_i^j$ are unknown coefficients to determine. In order to simplify the later description of the partitioned approach, we introduce the following notation for matrices and vectors: A function evaluated at a node set such as $\ut_j(X_j)$ denotes the column vector $(\ut_j(\x_1^j),\ldots,\ut(\x_{n_j}^j))^T$. An RBF evaluated at two node sets such as $\phi(Y_j,X_j)$ is an $(m_j\times n_j)$ matrix with elements $\phi(\|\underline{y}_i^j-\x_{k}^j\|))$, $i=1,\ldots,m_j$, $k=1,\ldots,n_j$, while $\phi(\x,X_j)$ is a row vector and $\phi(Y_j,\x)$ is a column vector. We can now write the RBF approximation~\eqref{eq:RBF} as
\begin{equation}
\tilde{u}_j(\x)={\phi}(\x,X_j)\Lambda_j,
\label{eq:RBF2}
\end{equation}
where ${\Lambda_j}=(\lambda_1^j,\ldots,\lambda_{n_j}^j)^T$.
It has been shown, e.g., in~\cite{DriFo02,LaFo05,Scha05}, that for infinitely smooth RBFs involving a shape parameter, the magnitude of the coefficients ${\Lambda_j}$ becomes unbounded as $\ep{}\rightarrow 0$ (flat limit), while the approximation $\tilde{u}_j(\x)$ itself remains well behaved. Therefore, we prefer to express the approximations in terms of the nodal values $\tilde{u}_j(X_j)$. We can use~\eqref{eq:RBF2} to form a linear system relating the coefficients to the nodal values
\begin{equation}
\phi(X_j,X_j){\Lambda_j}={\tilde{u}_j(X_j)}.
\label{eq:A}
\end{equation}  
For positive definite RBFs as well as for the multiquadric RBF, the interpolation matrix is non-singular for distinct node points~\cite{Scho38,Micchelli86}. By formally solving for $\Lambda_j$ in~\eqref{eq:A}, we can reformulate~\eqref{eq:RBF2} in terms of the nodal values as
\begin{equation}
\tilde{u}_j(\x)={\phi}(\x,X_j)\phi(X_j,X_j)^{-1}{\tilde{u}_j(X_j)}.
\label{eq:RBF3}
\end{equation}
This form also provides a definition of the nodal or cardinal basis as
$\psi(\x,X_j)={\phi}(\x,X_j)\phi(X_j,X_j)^{-1}$. Applying a linear differential operator to the RBF approximation results in
\begin{equation}
\mathcal{L}\tilde{u}_j(\x)=\mathcal{L}{\phi}(\x,X_j)\phi(X_j,X_j)^{-1}{\tilde{u}_j(X_j)}.
\label{eq:RBF4}
\end{equation}
Finally, we define a differentiation matrix
\begin{equation}
D^{\mathcal{L}}(Y_j,X_j) = \mathcal{L}\phi(Y_j,X_j)\phi(X_j,X_j)^{-1},
\label{eq:RBF5}
\end{equation}
such that
\begin{equation}
\mathcal{L}\tilde{u}_j(Y_j) = D^{\mathcal{L}}(Y_j,X_j)\tilde{u}_j(X_j).
\label{eq:RBF6}
\end{equation}
As mentioned above, we work in the nodal basis in order to avoid ill-conditioning as the shape parameter $\ep{}\rightarrow 0$ and the basis functions become increasingly flat. However, if we use~\eqref{eq:RBF4} directly to compute the differentiation matrices, we still need to deal with the ill-conditioning of $\phi(X_j,X_j)$. Instead, we use the RBF-QR approach~\cite{FoPi07,FoLaFly11,LLHF13}, which is a stable evaluation method that allows us to compute differentiation matrices for any small value of $\ep{}$. As shown in~\cite{LaHer17}, using RBF-QR or another method that is stable as $\ep{}\rightarrow0$~\cite{FoWri04,FassMcC12,FoLehPo13,WriFo17} is also vital in order to have convergence in an RBF partition of unity method.

\subsection{The RBF partition of unity method} In RBF-PUM, we combine the partition of unity approach with local RBF approximations. 
As we are aiming to solve a PDE problem, we now consider applying differential operators to the global partition of unity approximation~\eqref{eq:PUM} with the local RBF approximations~\eqref{eq:RBF3},
\begin{equation}
\mathcal{L}\ut(\x)=\sum_{j=1}^P\mathcal{L}\left(w_j(\x)\ut_j(\x)\right)
=\sum_{j=1}^P\mathcal{L}\left(w_j(\x)\phi(\x,X_j)\right)\phi(X_j,X_j)^{-1}\ut_j(X_j).
\label{eq:PUM2}
\end{equation} 
To fully expand the right hand side, we need to apply a product derivative rule. We make an example with the Laplacian, which is the operator we will use here
\begin{align}
\Delta\ut(\x)=\sum_{j=1}^P\left(\right.&\Delta w_j(\x)\phi(\x,X_j)+2\nabla w_j(\x)\cdot\nabla\phi(\x,X_j)\nonumber\\
&\left.+w_j(\x)\Delta \phi(\x,X_j)\right)\phi(X_j,X_j)^{-1}\ut_j(X_j),
\label{eq:PUM3}
\end{align} 
where the scalar product should be applied to the components of the gradients. To express this using differentiation matrices, we also need to put the weight function contributions into a proper matrix form. We let 
\[V_j^\mathcal{L}(Y)=\diag(\mathcal{L}w_j(Y)),\]
where $Y$ is an arbitrary set of evaluation points.
Then we can write the Laplacian as
\begin{equation}
\Delta\ut(Y)=\sum_{j=1}^P\left(V_j^\Delta(Y) D^I(Y,X_j) +2 V_j^{\nabla}(Y)\cdot D^{\nabla}(Y,X_j) + V_j^I(Y) D^\Delta(Y,X_j)\right)\ut_j.
\label{eq:PUM6}
\end{equation} 

\subsection{The two RBF-PUM formulations}
\begin{figure}
\centering
\includegraphics[width=0.42\textwidth]{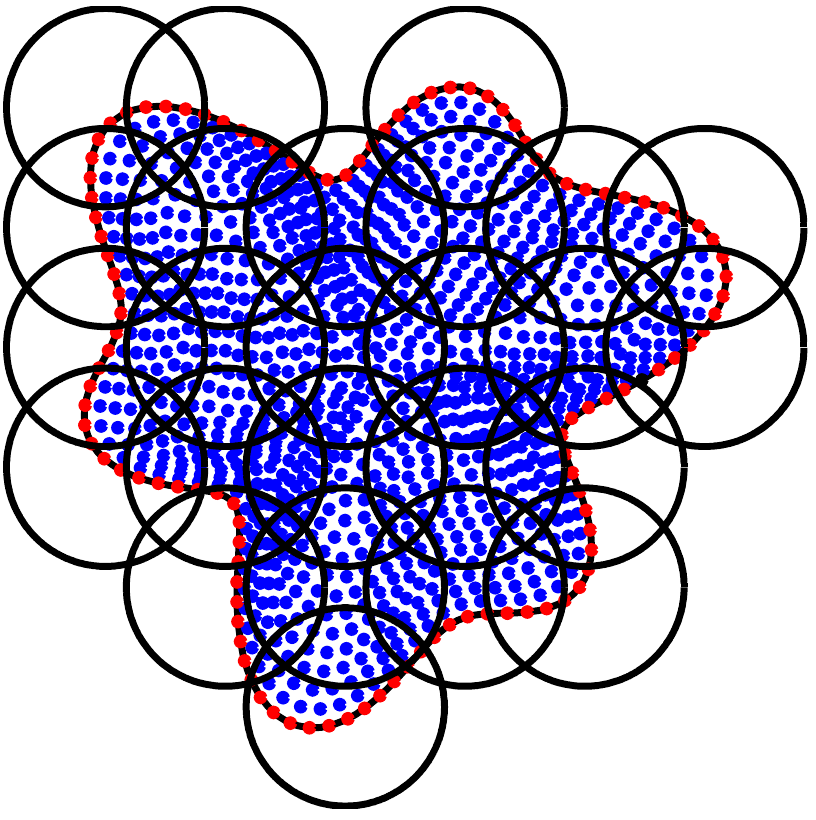}
\raisebox{5mm}{\includegraphics[width=0.34\textwidth]{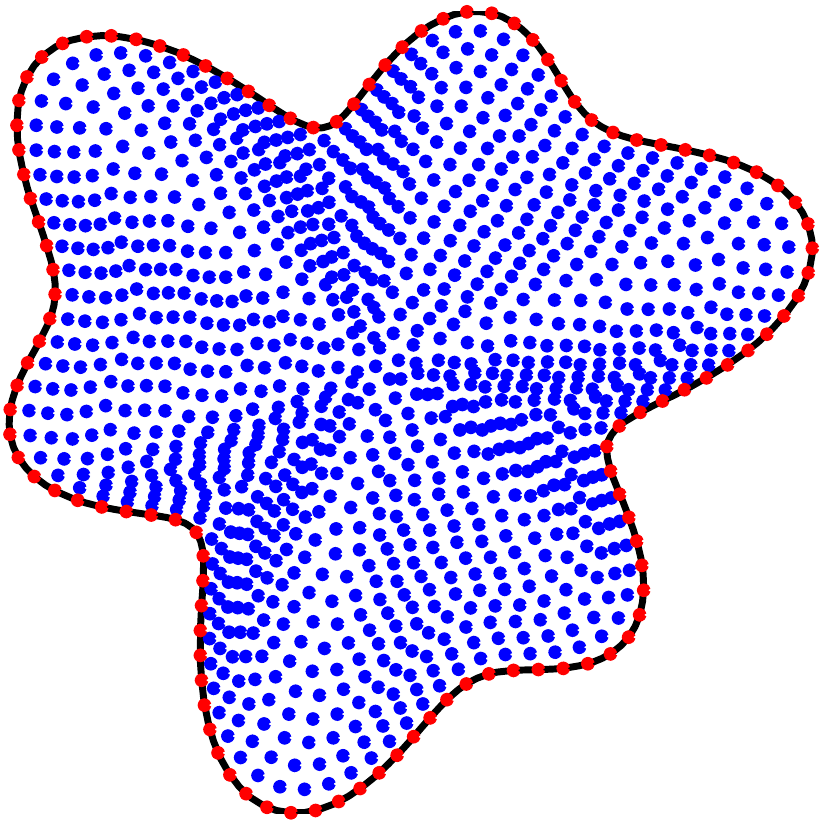}}
\caption{A global node set and patches in $\Omega_S$ (left), and evaluation points for collocation (same as node points) (right) for C-RBF-PUM.}
\label{fig:collnodes}
\end{figure}

The main differences between C-RBF-PUM and LS-RBF-PUM can be explained through the choice of node points and evaluation points. In the C-RBF-PUM formulation, a global node set $X=\{\x_k\}_{k=1}^{N}$ with a subset of interior points $X^i=\{\x_k\in X : \x_k\in\Omega\}$ and of boundary points $X^b=\{\x_k\in X : \x_k\in\partial\Omega\}$ is constructed. The local node sets $X_j$ are then derived from the global set. The evaluation (collocation) point set $Y=X$, see Figure~\ref{fig:collnodes}.

The number of points $n_j$ in the local node sets $X_j=Y_j$ is similar for interior patches, but can vary significantly for patches that intersect the boundary. Using the relative area, for a quasi uniform node distribution we have
\[n_j\approx N\frac{|\Omega_j\cap\Omega|}{|\Omega|}.\]
That is, patches with a small intersection can have a much lower number of local points, which results in a lower approximation order. As is further discussed in Section~\ref{sec:theory} this reduces the global convergence rate.

\begin{figure}
\centering
\includegraphics[width=0.42\textwidth]{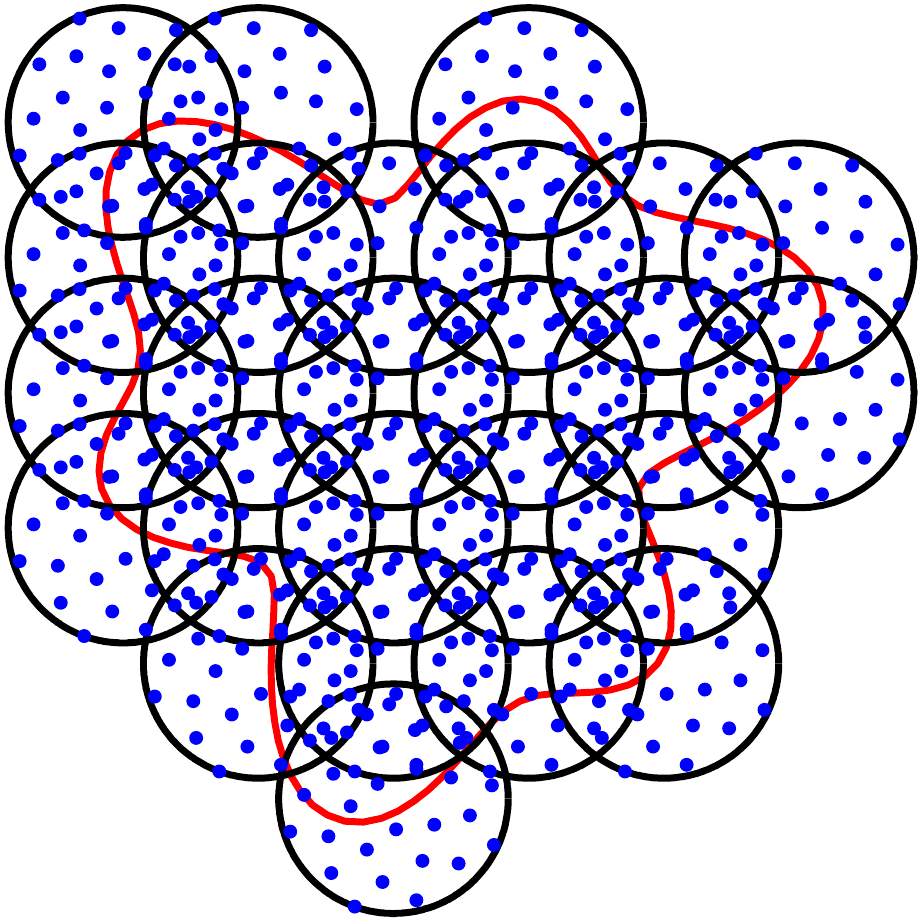}
\raisebox{5mm}{
\includegraphics[width=0.34\textwidth]{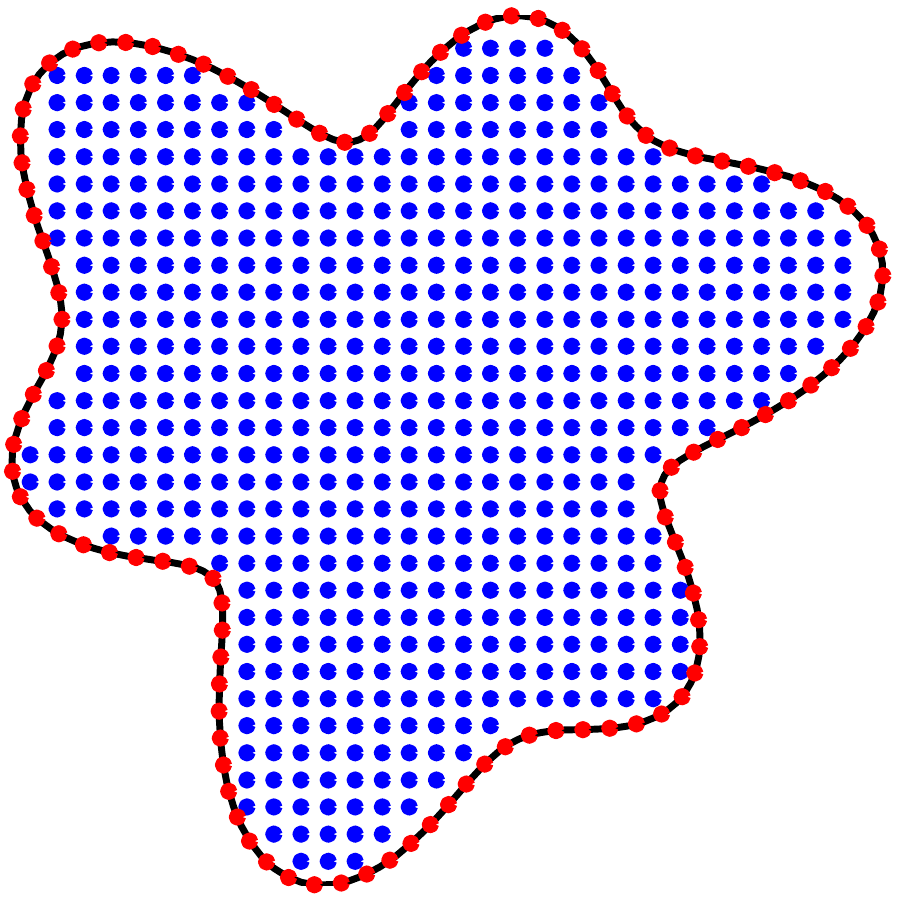}
}
\caption{Patches with identically distributed local node sets covering the domain $\Omega_S$ (left), and least squares evaluation points on a Cartesian grid in the interior and uniform with respect to arc length on the boundary (right) for LS-RBF-PUM.}
\label{fig:nodes}
\end{figure}

For LS-RBF-PUM, we instead start from local node sets $X_j$ that are identically distributed with respect to the corresponding patches $\Omega_j$. This allows us to create optimized node sets for discs and spheres that can be used for any geometry of the domain $\Omega$. The global node set is here the union of the local sets, $X=\bigcup_{i=1}^PX_j$. For each patch, the number of local points $n_j= n$, and the global number of node points is $N=nP$. We completely decouple the least squares evaluation points $Y$ from the node points. This allows us to choose a simple scheme for the layout. We use evaluation points $Y^i\subset Y$ distributed on a Cartesian grid inside the domain and points $Y^b\subset Y$ distributed uniformly with respect to the arc length on the boundary of the domain. An example of LS-RBF-PUM node and patch layouts is shown in Figure~\ref{fig:nodes}.

Some advantages of LS-RBF-PUM over C-RBF-PUM follow directly from the choice of nodes.
\begin{itemize}
\item 
Having the same number of local points in all patches ensures a similar approximation order in the whole domain, see also Section~\ref{sec:theory}.
\item 
High quality node points can be generated independently of the problem geometry and a simple scheme can be used for the evaluation points.
\item 
Referring, e.g., to~\eqref{eq:PUM2}, we see that the inverse of the local interpolation matrix $\phi(X_j,X_j)$ is needed for each patch. Especially when using RBF-QR, forming and factorizing these matrices is a costly operation. When we are using identical distributions in LS-RBF-PUM, we only need to do this for one patch, which significantly reduces the setup cost.
\end{itemize}
Note that we allow node points to fall outside of the domain. The corresponding basis functions contribute to the solution inside the domain. Also note that for LS-RBF-PUM the local solutions in two overlapping patches do not need to match in the overlap region, while for C-RBF-PUM, we enforce unique values at collocation points in overlap regions.

\subsection{The method parameters} 
When deriving theoretical estimates and performing numerical experiments, we use a number of key parameters to describe the method. For defining the patches, we use an underlying box structure, and let each patch be centered in the box center. The side length of the boxes is denoted by $H$. Patches that circumscribe their respective boxes have no overlap in the diagonal direction. This corresponds to a radius $\rho_0=\sqrt{d}H/2$. The overlap parameter $\delta$ describes the overlap between patches. The radii of patches with overlap $\delta$ are given by $\rho=(1+\delta)\rho_0$. The number of patches $P$ is determined by the choice of $H$ and $\delta$, through the intersection of the generated patches with the domain. 

A node set is characterized by its fill distance $h$, which corresponds to the radius of the largest ball fully contained in $\Omega$ that is empty of node points, 
\begin{equation}
h=\sup_{\x\in\Omega}\ \min_{\x_j\in X}\|\x-\x_j\|.
\label{eq:filldist}
\end{equation}
Together with the choice of distribution scheme for the nodes, $h$ determines the number of local points $n_j$ within each patch. For LS-RBF-PUM, $n_j\equiv n$. For C-RBF-PUM, we let this parameter measure the worst case over all patches, $n=\min_jn_j$. 

For LS-RBF-PUM, the rate of oversampling $\beta=M/N$, where $M>N$ is the number of least squares evaluation points and $N$ is the number of node points, is an important parameter for the performance of the method. 

Finally, the shape parameter $\ep{}$ of the RBFs, that govern their relative flatness, is relevant for the accuracy of the approximations.




\subsection{Solving the Poisson problem using RBF-PUM}
\begin{figure}
\centering
\includegraphics[height=0.4\textwidth]{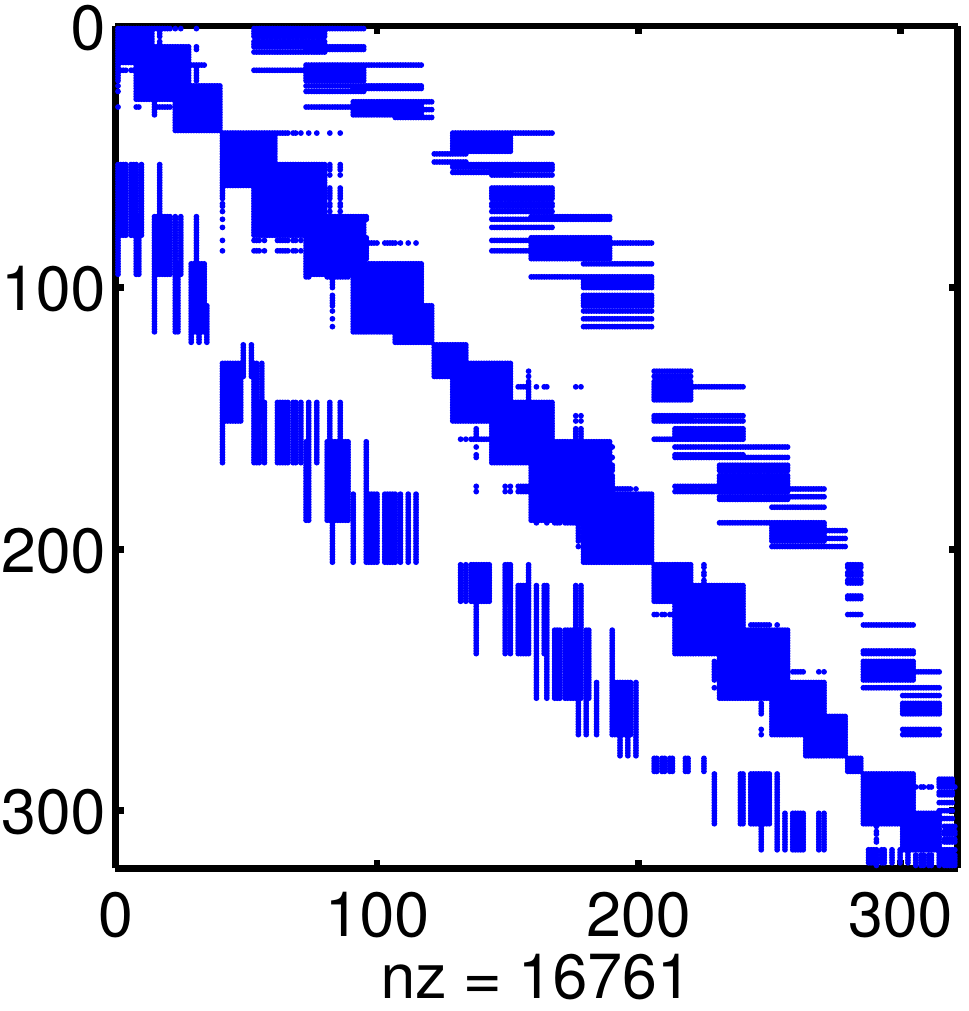}
\includegraphics[height=0.4\textwidth]{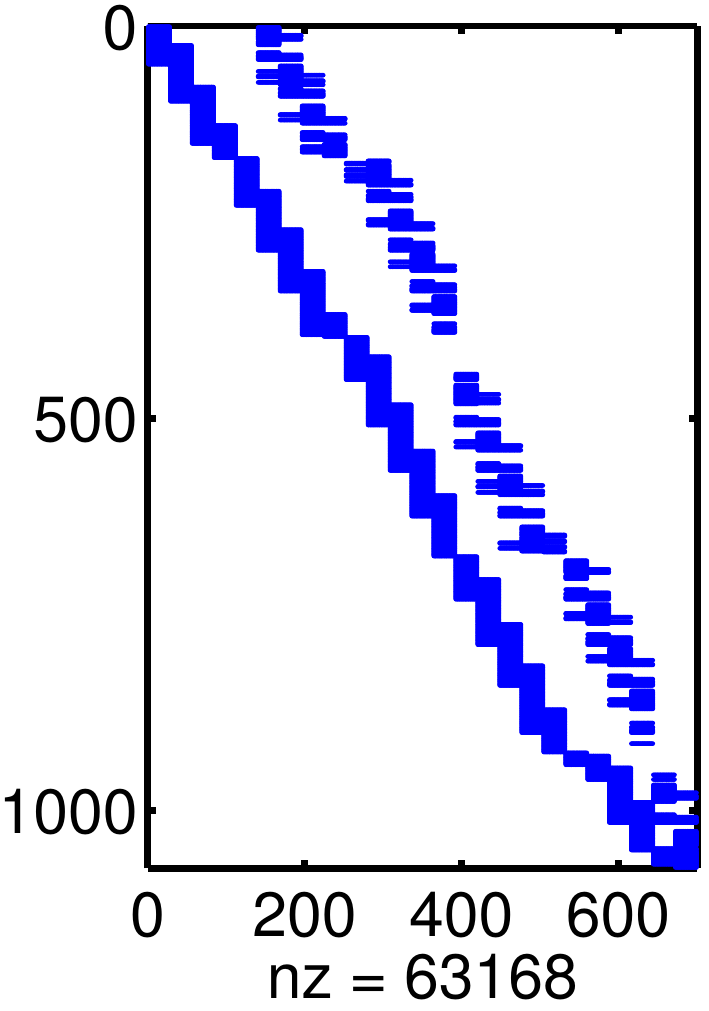}
\caption{The structure of the matrix $L$ for a problem defined over $\Omega_S$ with $P=24$ patches, box size $H=0.6$, and overlap $\delta=0.2$. For C-RBF-PUM (left), $13\leq n_j\leq42$, $h\approx 0.12$, and $N=321$, and for LS-RBF-PUM (right), $n=28$, $h\approx 0.14$, $N=700$, $M=1073$, and $\beta=M/N\approx 1.5$.}
\label{fig:struct}
\end{figure}

To solve~\eqref{eq:poisson}, we set up a linear system, where each evaluation point (test point) corresponds to one equation. We enforce the boundary condition for points in $Y^b$, and we enforce the PDE for points in $Y^i$. 
We assemble the global matrix by adding the contributions from each patch. For patch $j$, we compute the following matrix blocks:

\[L_j^i = V_j^\Delta(Y_j^i) D^I(Y_j^i,X_j) +2 V_j^{\nabla}(Y_j^i)\cdot D^{\nabla}(Y_ j^i,X_j) + V_j^I(Y_j^i) D^\Delta(Y_j^i,X_j),\]
\[L_j^b = V_j^I(Y_j^b) D^I(Y_j^b,X_j).\]
The global matrix $L$ has size $(M\times N)$ in the case of LS-RBF-PUM and $(N\times N)$ for C-RBF-PUM. To add the local contributions to the global matrix we need the indices $I_j^i$ and $I_j^b$ of the local evaluation points  in the global evaluation point set $Y$ as well as the indices $J_j$ of the local node points in the set $X$. These form the indices in the global matrix, which we assemble as
\begin{align}
L(I_j^i,J_j) &= L(I_j^i,J_j) + L_j^i, \quad j=1,\ldots,P,\\
L(I_j^b,J_j) &= L(I_j^b,J_j) + L_j^b, \quad j=1,\ldots,P.
\end{align}
The global matrix is sparse, which is crucial to scale the method to large problem sizes without prohibitive computational cost and memory requirements. The structure of the matrix depends on the ordering of the nodes. In~\cite{HLRvS16}, a particular ordering was used to provide a structure suitable for preconditioning. Here, we order the nodes according to patch in a greedy sense such that we start with $Y_1$, then we take all points in $Y_2$ that were not already picked ($Y_2\setminus Y_2\cap Y_1$) until we pick the last points left in $Y_P$. The type of structures that arise are illustrated in Figure~\ref{fig:struct}.

The global linear system to solve is
\begin{equation}
L(Y,X)U(X)=F(Y),
\end{equation}
where $U(X)=\ut(X)$ for C-RBF-PUM and
\[U(X)=\left(\begin{array}{c}
 \ut_1(X_1)\\
\vdots\\
\ut_P(X_P)
\end{array}\right)\]
for LS-RBF-PUM, and
\[
F(\y)=\left\{\begin{array}{ll}
f(\y), &\y\in \Omega,\\
g(\y), & \y\in \partial\Omega.
\end{array}\right.\]
For C-RBF-PUM, $Y=X$, and we have a square linear system to solve either through direct factorization or through an iterative method~\cite{HLRvS16}. In this paper we use LU-factorization to provide a fair comparison.

 For LS-RBF-PUM, we solve the system using QR-factorization. To implement an iterative approach is a potential future development. We factorize the matrix 
\begin{equation}
L(Y,X)=\left(\begin{array}{cc}
\\
Q_1 & Q_2\\
\\
\end{array}\right)
\left(\begin{array}{c}
R_1\\0
\end{array}\right),
\end{equation}
where the $(M\times N)$ matrix $Q_1$ forms an orthogonal basis for the span of the columns of $L(Y,X)$, and the $(N\times N)$ matrix $R_1$ is upper triangular. The least squares solution is obtained from solving the system
\begin{equation}
R_1U(X)=Q_1^TF(Y).
\end{equation}
The residual can be expressed as $r_U(Y)=L(Y,X)U(X)-F(Y)=Q_2Q_2^TF(Y)$. We note for later use that the following orthogonality relation holds
\begin{equation}
L^Tr_U=0.
\label{eq:ortho}
\end{equation}

We expect the number of least squares evaluation points $M$ to be larger than the number of node points $N$, but there is also one further requirement. The number of evaluation points $Y_j$ within a patch $\Omega_j$ must be larger than or equal to the number of node points $X_j$ in the same patch for $R_1$ to have full rank. This can be a problem for boundary patches that only contain a small part of $\Omega$, and thereby a relatively small ratio of test points (inside $\Omega$) to node points (whole patch). This issue can be resolved by refining, shifting and/or scaling boundary patches. If the node layout relative to the patch remains unchanged, these modifications can be incorporated without significantly increasing the computational cost of the method.

\section{Theory}
\label{sec:theory}
We start by defining the RBF-PUM interpolant which is used as an auxiliary function in the error estimates,
\begin{equation}
  \mathcal{I}(u) = \sum_{j=1}^Pw_j\mathcal{I}(u_j),
\label{eq:interp}
\end{equation}
where $\mathcal{I}(u_j)$ is the local RBF interpolant of form~\eqref{eq:RBF} satisfying $\mathcal{I}(u_j)(X_j)=u_j(X_j)$. We need estimates for the interpolation error and its derivatives. We define
\begin{equation}
\mathcal{E}_\mathcal{L} = \mathcal{L}(\mathcal{I}(u)-u).
\end{equation}
RBF-PUM interpolation errors were studied extensively in~\cite{LaHer17} using sampling inequalities from~\cite{RieZwi10}. Two types of estimates were provided.
\begin{estimate}[\cite{LaHer17,SVHL15}]
For an RBF-PUM interpolant to a function $u$ over a domain $\Omega$, where $n_j$ is kept fixed while the patch sizes $H_j$ are refined, the following algebraic error estimate holds:
%
%
\begin{equation}
\|\mathcal{E}_\mathcal{L}\|_{L_\infty(\Omega)}\leq K\max_{1\leq j\leq P}C^A_{j}H_j^{q(n_j)+1-\frac{d}{2}-\alpha}\|u\|_{\mathcal{N}({\Omega}_j)},
\label{eq:interperr1}
\end{equation}
where the constants $C^A_j$ depend on the dimension $d$, the chosen weight function, the number of local points $n_j$, and the order $\alpha$ of the differential operator. 
%
The function $q(n_j)$ corresponds to the polynomial degree $q$ supported by the local number of points $n_j$. Let $n_{q,d}$ be the dimension of the polynomial space of degree $q$ in $d$ dimensions.
If the number of local points satisfies $n_{q,d}\leq n_j < n_{q+1,d}$, then $q(n_j)=q$. The norm $\|\cdot\|_{\mathcal{N}({\Omega}_j)}$ is the native space norm~\cite{Fass07} of the space generated by the chosen RBFs.
\label{est:alg}
\end{estimate}
We do not go into details regarding unisolvency of the node sets here, as we have the option to choose appropriate nodes. It should also be noted that our test functions are not chosen to lie in the native space of the Gaussian RBFs. In practice approximation of smooth functions works well, and as the patch size is refined, both the local native space and the local function space approach a polynomial space. For further discussion of these topics, see~\cite{LaFo05}.

 In the estimate, we include the variation over the patches. We expect that adaptive approaches based on this and the following estimate will be of interest as a future development.
 
\begin{estimate}[\cite{LaHer17,SVHL15}]
For an RBF-PUM interpolant to a function $u$ over a domain $\Omega$, using Gaussian RBFs, where the patch size $H_j$ is kept fixed, while the local node density $h_j$ is varied, the following exponential error estimate holds when $h_j$ is sufficiently small:

\begin{equation}
\|\mathcal{E}_\mathcal{L}\|_{L_\infty(\Omega)}\leq KC^E\max_{1\leq j\leq P} e^{\gamma\log(h_j)/\sqrt{h_j}}\|u\|_{\mathcal{N}({\Omega}_j)}.
\label{eq:interperr2}
\end{equation}
where the constant $C^E$ and the rate $\gamma$ both depend on the dimension $d$ and the order of the differential operator, and $C^E$ additionally depends on the chosen weight function.
\label{est:spec}
\end{estimate}
This estimate is for Gaussian RBFs, but similar estimates can be constructed, e.g., for inverse multiquadrics~\cite{RieZwi10}. Here, we have not gone into details concerning how the patches $\Omega_j$ intersect the domain $\Omega$. For C-RBF-PUM, this affects the interior cone condition, which in turn affects the constants in the estimates. 

The interpolation error estimates are essentially the same whether we use a collocation approach or a least squares approach. The interpolation error drives the convergence, and hence, we cannot expect a different order of convergence due to the introduction of least squares testing. However, we do expect that the numerical robustness of the method for large problem sizes will be improved. In order to show this, we need to look at the full error estimate.

We start from the well-posedness, Estimate~\ref{est:well}, of the elliptic PDE. Then we insert the RBF-PUM interpolant~\eqref{eq:interp} as an auxiliary function, and finally we use the interpolation error estimate~\eqref{eq:interperr1} or~\eqref{eq:interperr2} to get
\begin{align}
\|\tilde{u}-u\|_{L_2(\Omega)}&\leq C_P \|\tilde{u}-u\|_F\nonumber\\
& \leq C_P\left( \|\tilde{u}-\mathcal{I}(u)\|_F+ \|\mathcal{I}(u)-u\|_F\right)\nonumber\\
& \leq C_P\left( \|\tilde{u}-\mathcal{I}(u)\|_F+ \|\mathcal{E}_\mathcal{L}\|_{L_\infty(\Omega)}\right).
\label{eq:est1}
\end{align}
The term that remains to be estimated contains the RBF-PUM approximant and the RBF-PUM interpolant. Therefore, the operator implied by the data norm $\|\cdot\|_F$ can be applied through a common differentiation matrix. We have
\begin{equation}
\mathcal{L}(\tilde{u}-\mathcal{I}(u)) = L(\x,X)(U(X)-u(X)).
\label{eq:cont}
\end{equation}
However, we need to transform this further to produce an estimate. We have information about the residual at the evaluation set $Y$. We start by noting that
\begin{equation}
\mathcal{L}(\tilde{u}(Y)-\mathcal{I}(u)(Y)) =  L(Y,X)(U(X)-u(X)).
\label{eq:disc}
\end{equation}
Assuming that $L\equiv L(Y,X)$ has full column rank (this is expected), we can construct a pseudo inverse $L^+=(L^TL)^{-1}L^T$ such that $L^+L=I$.

Inserting $L^+L$ into~\eqref{eq:cont} and using~\eqref{eq:disc} we can rewrite~\eqref{eq:cont} as
\begin{align}
\mathcal{L}(\tilde{u}-\mathcal{I}(u)) &=  L(\x,X)L^+\mathcal{L}(\tilde{u}(Y)-\mathcal{I}(u)(Y))\nonumber \\
& = L(\x,X)L^+\left(\mathcal{L}(\tilde{u}(Y)-u(Y))-\mathcal{L}(\mathcal{I}(u)(Y)-u(Y))\right)\nonumber \\
& = L(\x,X)L^+\left(r_U-\mathcal{E}_\mathcal{L}(Y)\right).
\label{eq:est2}
\end{align}
We now use the orthogonality property~\eqref{eq:ortho} of the residual, which means that in exact arithmetic $L^+r_U=0$. We are however interested in the effects of numerical errors. Therefore, we assume that the orthogonality relation holds to within a multiple $C_M$ of the machine precision $\delta_M$. In the case of C-RBF-PUM, the linear system is square and the residual itself should be on the order of $\delta_M$.

Combining~\eqref{eq:est1}, \eqref{eq:est2}, and replacing the residual term with the rounding error, we get the final estimate   
\begin{equation}
\|\tilde{u}-u\|_{L_2(\Omega)}\leq C_P\|\mathcal{E}_\mathcal{L}\|_{L_\infty(\Omega)} + C_P
\|L(\cdot,X)L^{+}\|_{L_\infty(\Omega)}\left( C_M\delta_M + \|\mathcal{E}_\mathcal{L}\|_{L_\infty(\Omega)}\right).
\label{eq:errest}
\end{equation}
The error estimate is proportional to the interpolation error as expected, down to the lower limit provided by the rounding error. The matrix norm involving $L$ is similar to a condition number for the PDE approximation. How it correlates with the problem size and problem parameters is important for the numerical robustness of the method.


We have not managed to provide a theoretical bound for the matrix norm that allows for convergence. However, we know from~\cite{Schaback16} that when a nodal basis is used and oversampling is employed on the test side, uniform stability (no growth of the matrix norm) can be achieved.
An important issue to quantify is how much oversampling is needed.
In the following section, we investigate the matrix norm, which we call the stability norm, numerically.
\section{Numerical experiments}\label{sec:exp}
Both of the RBF-PUM algorithms are implemented in MATLAB. The numerical experiments for the two-dimensional test cases are carried out on a MacBook Pro with Core i7 and 16 GB RAM, while the experiments for the three-dimensional test cases are performed at the UMass Dartmouth rapid prototyping server, a dual 8-core Intel Xeon 2.2 GHz workstation with 32 GB RAM. The most extensive tests are carried out for the two-dimensional problems, and then in the final subsection we verify that the method behaves as expected also in the three-dimensional case.

The RBF-QR method for stable evaluation of the local differentiation matrices~\cite{FoLaFly11,LLHF13} is used in all experiments. The overlap parameter is set to $\delta=0.2$, a choice which in our experiments has shown to be effective. Using a smaller overlap parameter increases the error, while a larger overlap increases the amount of work. For $u_1$, $u_2$ and $u_5$, the shape parameter $\ep{}=1$ is used, and for $u_3$ and $u_4$, $\ep{}=4$. The default value for the rate of oversampling $\beta\approx1.5$ for all experiments. Unless otherwise stated, the experiments in two dimensions are performed on the square computational domain $\Omega_B=[-2,\,2]^2$. By choosing a regular domain as the square for analysis of the method performance, we eliminate noise due to variability of the geometry in relation to the patch layout. Irregular domains are investigated in a separate subsection. 

For C-RBF-PUM, we use a uniform Cartesian node distribution. This is a good choice from the point of view that the nodes are uniform and easy to generate. However, it also leads to sensitivity regarding the alignment with the patches as can be seen in the experiments. 

\begin{figure}[!htb]
\centering
\includegraphics[width=0.2\textwidth]{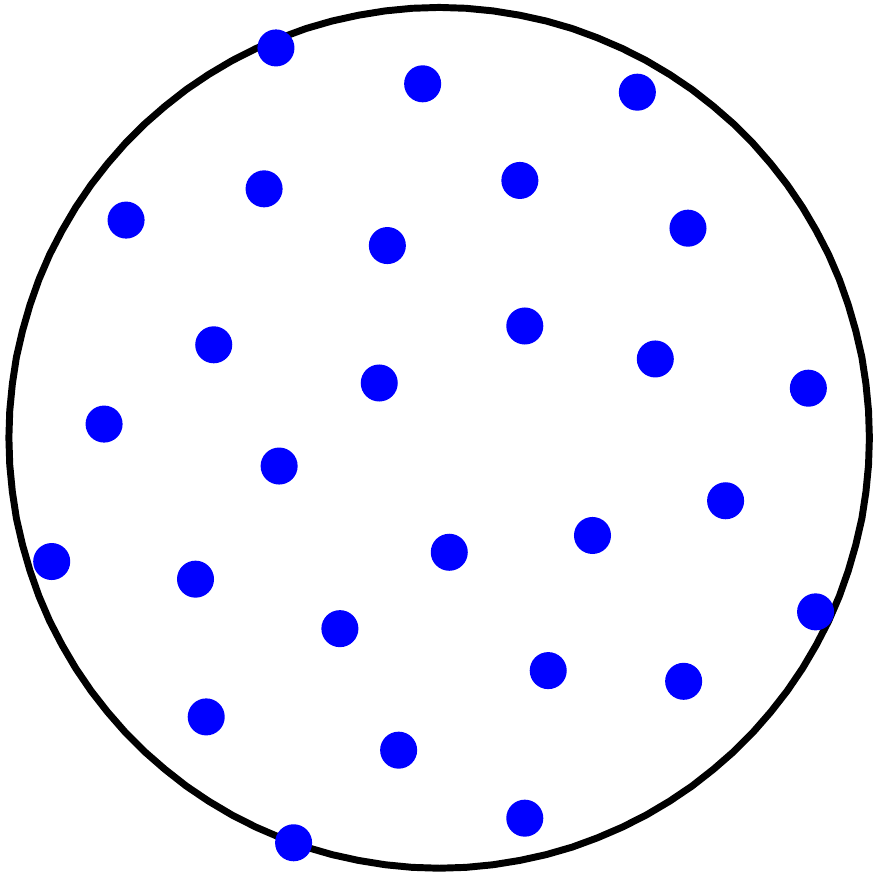}\hspace*{0.05\textwidth}
\includegraphics[width=0.2\textwidth]{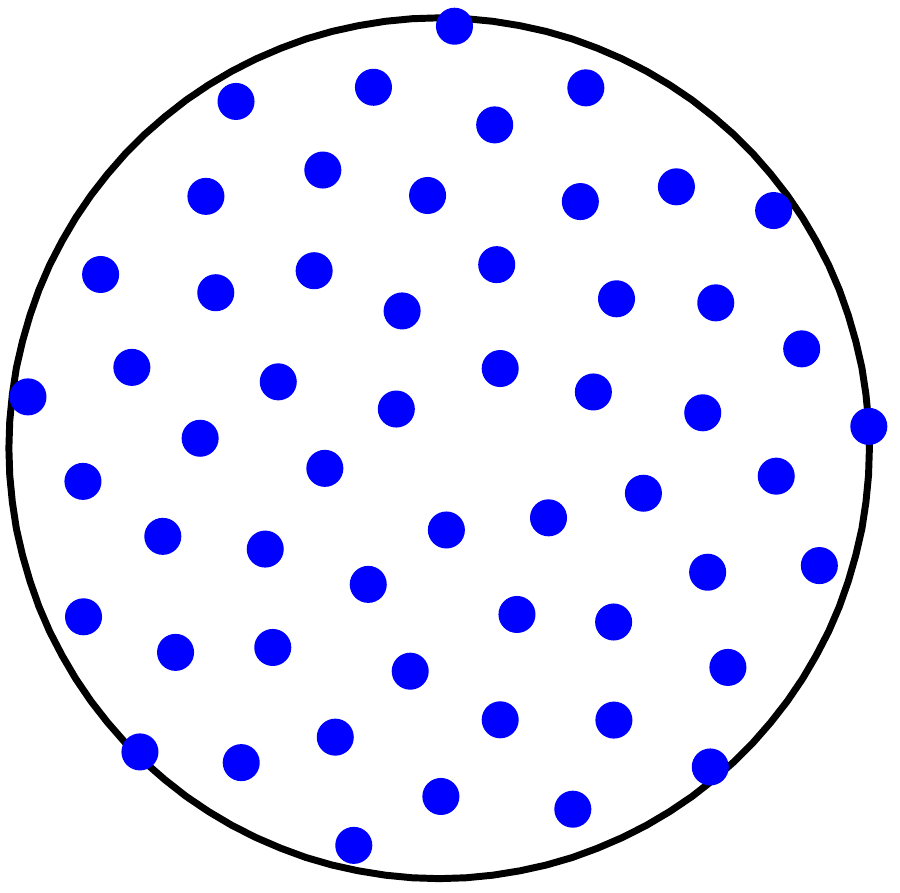}\hspace*{0.05\textwidth}
\includegraphics[width=0.2\textwidth]{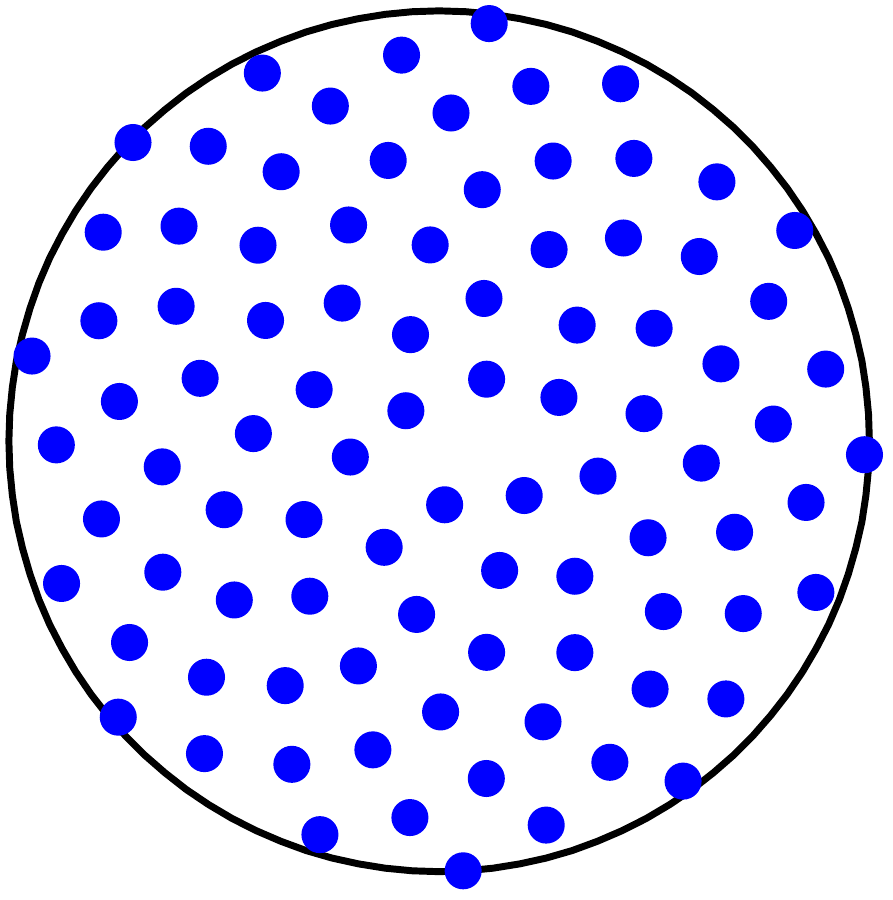}
\caption{The layout of $n=28$, $n=55$ and $n=91$ Vogel points in a patch.}
\label{fig:vogel}
\end{figure}

For LS-RBF-PUM, in the two-dimensional case, we use a Vogel node distribution $\x_i=\sqrt{i/n}\left(\cos(i\hat{\theta}),\sin(i\hat{\theta})\right)$, $i=1,\ldots,n$, where $\hat{\theta}=\pi(3-\sqrt{5})$, in each patch. These nodes are quasi uniform and we can easily control $n$, see Figure~\ref{fig:vogel} for some examples. We have also tried other types of node sets in the disc, including nodes clustered toward the boundary, but we did not observe any significant differences in the results.

For the three-dimensional case, we use 
an ad hoc optimization approach to generate a quasi uniform mesh in a sphere and then we use these nodes for each patch. For each new point, we minimize the distance to the origin, under the constraint that the distance to the nearest neighbor is $\geq1$. The resulting node set is then scaled to fit the patch size. A spherical patch with a local node set is illustrated in Figure~\ref{fig:dom3}. 

Examples of evaluation points in two and three dimensions are provided in Figures~\ref{fig:nodes} and~\ref{fig:dom3}. The evaluation points are distributed on a Cartesian grid in the interior of $\Omega$ and then uniformly or quasi uniformly on $\partial\Omega$.

The theoretical results are derived in the $L_2$-norm, but we measure the errors in the $L_\infty$-norm, since this implies that the result holds also in the $L_2$-norm with an additional constant. The error as well as the stability norm estimates are evaluated at 1000 Halton nodes. This error measure provides an estimate from below of the continuous norm, but we have verified that the number of points is large enough that the difference is small compared with a more dense sampling.


\subsection{Numerical convergence results for the two-dimensional Poisson problem} The convergence of both the collocation and least squares formulation of RBF-PUM is governed by the interpolation error. The aim of these experiments is to see if the numerical convergence behavior follows the theoretical predictions as well as to compare the results for the two formulations. 

\begin{figure}[!htb]
\centering
\includegraphics[height=0.38\textwidth]{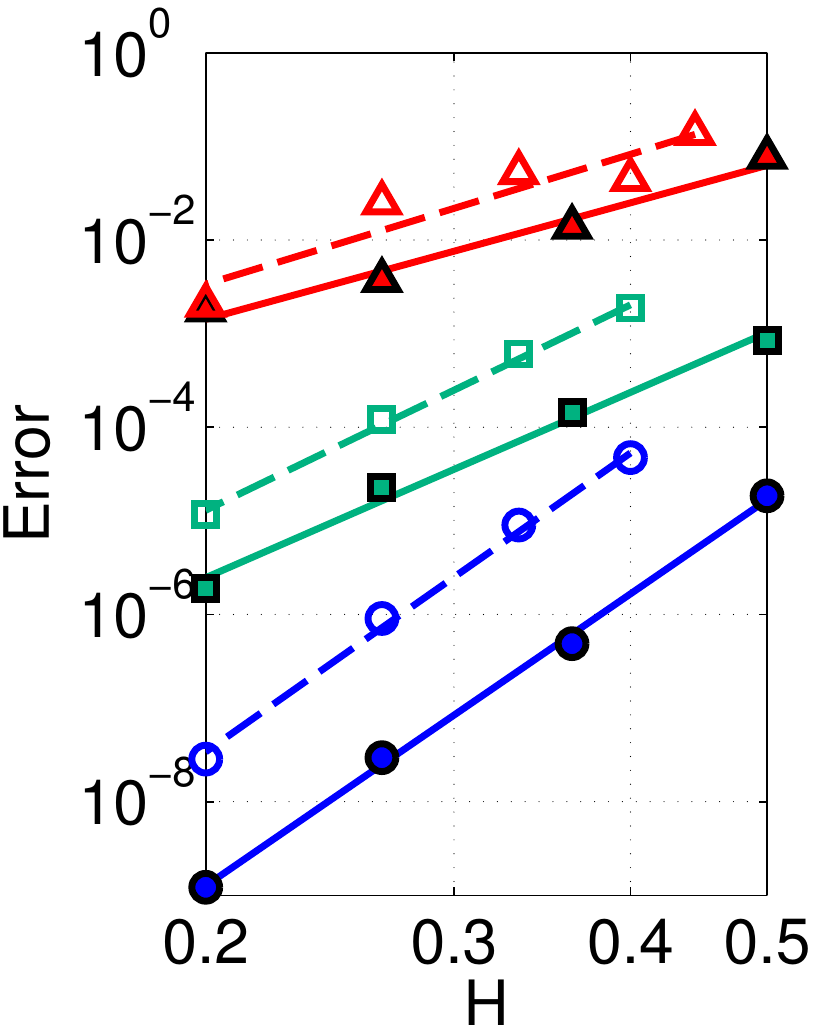}
\includegraphics[height=0.38\textwidth]{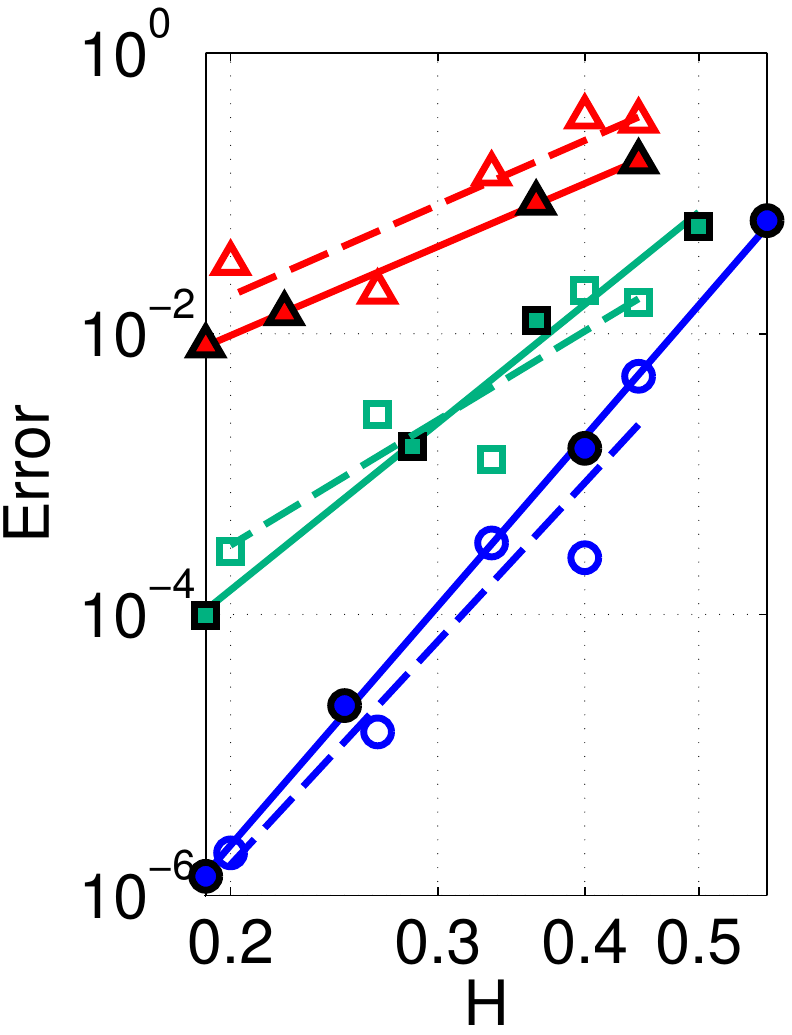}
\caption{Algebraic convergence of the error as a function of the patch size $H$ for fixed numbers of points per patch $n=28$ ($\triangle$), $n=55$ ($\Box$), and $n=91$ ($\bigcirc$), for collocation (dashed lines, open markers) and least squares (solid lines, solid markers) for the trigonometric function $u_2$ (left) and the Runge function $u_3$ (right). The numerically estimated slopes are $p=4.6,\,7.3,\,11.0$ for C-RBF-PUM for $u_2$, $p=4.1,\,6.6,\,10.0$ for LS-RBF-PUM for $u_2$, $p=3.7,\,5.1,\,9.0$ for C-RBF-PUM for $u_3$, and $p=3.6,\,6.8,\,9.7$ for LS-RBF-PUM for $u_3$.}
\label{fig:algconv}
\end{figure}

In the first experiment, we fix the number of points per patch $n$. Note that for C-RBF-PUM, due to the Cartesian node layout and boundary effects, this can only be done approximately. 
Theory predicts algebraic convergence in the patch size, see~\eqref{eq:interperr1}.
The numbers of local points in the experiment are chosen to be $n=28$, 55, 91, corresponding to convergence orders $p=4$, 7, and 10 in Estimate~\ref{est:alg}. Figure~\ref{fig:algconv} shows the maximum error as a function of the patch size. The slopes estimated from the numerical results through linear regression for a certain $n$ are similar for both methods and both test functions. The average numerically estimated rates are $\tilde{p}=4.0$, 6.4, and 9.9. That is, they are very close to the theoretical results.
The accuracy of the two methods is similar, but the results are more irregular for C-RBF-PUM. Furthermore, here only the best combinations of $h$ and $H$ are used for C-RBF-PUM in order to observe a convergence trend. The results for C-RBF-PUM are closer to those of LS-RBF-PUM for $u_3$ than for $u_2$. This may be explained by the fact that the Runge function $u_3$ is small near the boundary, where the C-RBF-PUM approximation may be less accurate due to the intersection of patches with the boundary.

\begin{figure}[!htb]
\centering
\includegraphics[height=0.38\textwidth]{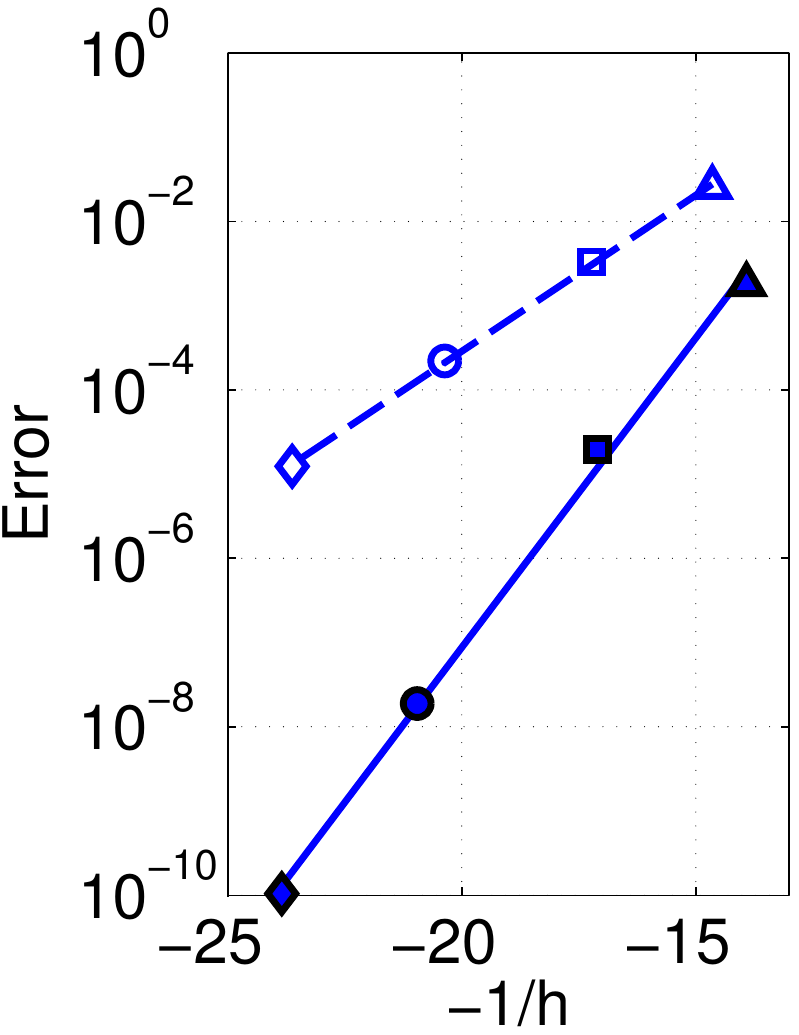}
\includegraphics[height=0.38\textwidth]{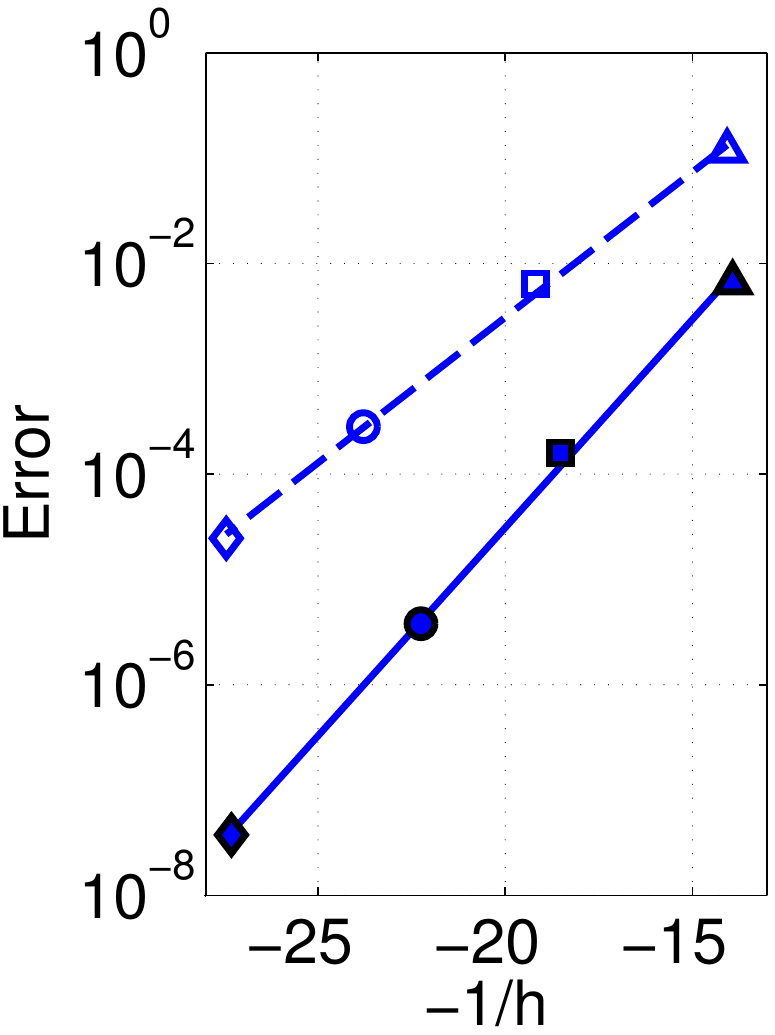}
\caption{Spectral convergence of the error as a function of the (negative) inverse fill distance $-h^{-1}$ for fixed patch size $H=0.2$ and $n=28$, $55$, $91$, $153$ for collocation (dashed lines, open markers) and least squares (solid lines, solid markers) for the trigonometric test case $u_2$ (left) and the Runge function $u_3$ (right).}
\label{fig:expconv}
\end{figure}

In the second experiment, we fix $H=0.2$ resulting in a total of $P=400$ patches, and then let $n$ vary. The values that are used are $n=28$, 55, 91, and 153, resulting in $N=11\,200$, 22\,000, 36\,400, and 61\,200 nodes for LS-RBF-PUM. 
Figure~\ref{fig:expconv} shows the corresponding convergence results. The horizontal axis corresponds to the inverse fill distance in order to illustrate spectral convergence of the form $\exp(-\gamma/h)$. This is not exactly the form in Estimate~\ref{est:spec}, but this is the behavior that we observe in practice. The fit to a line is equally good with $\log(h)/h$, but possibly a little worse with $\log(h)/\sqrt{h}$. Here, the accuracy of LS-RBF-PUM is significantly higher than for C-RBF-PUM. The main reason is that for a given fill distance $h$, the number of nodes per patch in the worst case for C-RBF-PUM is much lower in the corner and boundary patches. Therefore, the global accuracy is reduced compared with LS-RBF-PUM, where all patches have the same number of nodes.
%

Going back to the overall error estimate~\eqref{eq:errest}, we see that another candidate for differences in behavior between the methods is the generalized matrix norm associated with bounding the continuous operator in terms of the discrete operator. This norm is investigated numerically in the following subsection.


\subsection{Approximation stability}
\label{sec:stab}
 
\begin{figure}[!htb]
\centering
\includegraphics[height=0.36\textwidth]{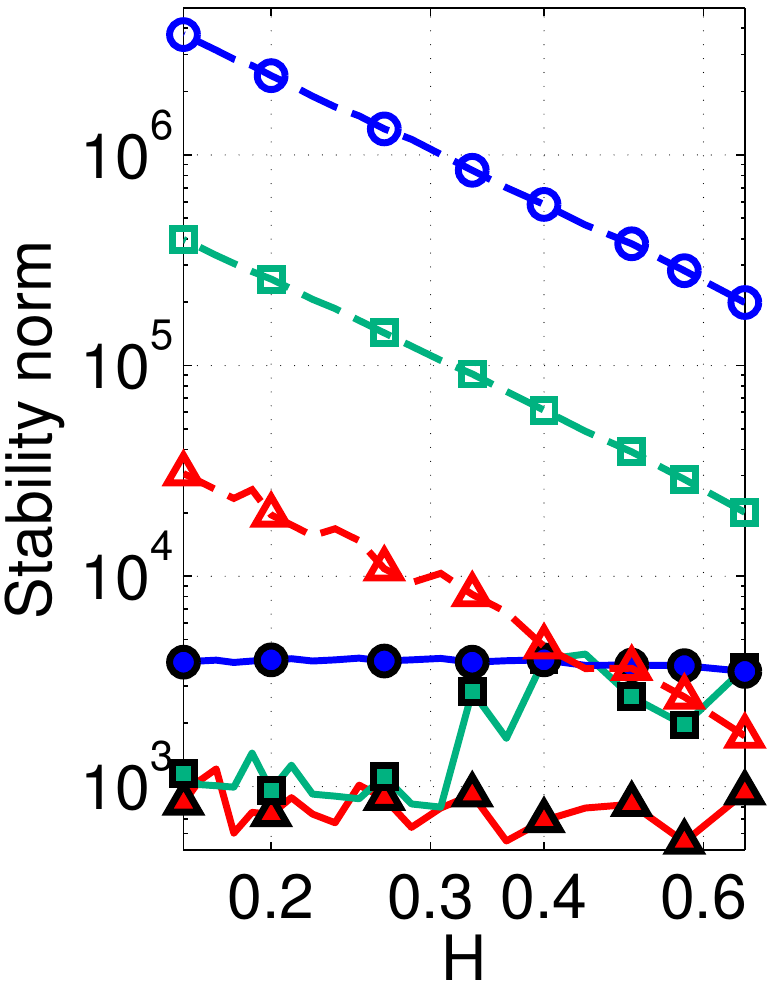}
\includegraphics[height=0.36\textwidth]{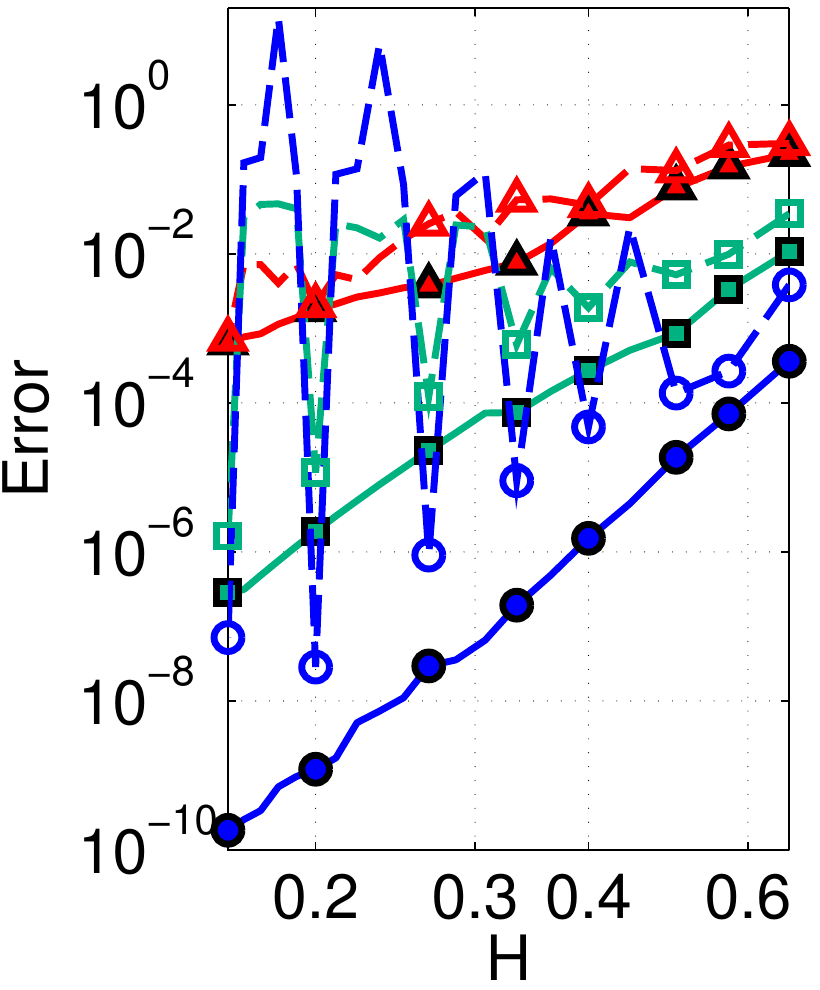}
\caption{The numerically estimated stability norm (left) and the corresponding error (right) as a function of $H$ for $n=28$ ($\triangle$), $n=55$ ($\Box$), and $n=91$ ($\bigcirc$) for LS-RBF-PUM (solid lines, solid markers) and C-RBF-PUM (dashed lines, open markers). Note that only selected data points, those that are optimal for C-RBF-PUM, have markers.}
\label{fig:boxnorm}
\end{figure}

We start the investigation of the stability norm from the case with a fixed number of nodes per patch and varying patch size. 
The resulting norm estimates as well as the errors in the solution are shown in Figure~\ref{fig:boxnorm}. When a collocation approach is used, the stability norm grows algebraically as the patch size $H$ decreases. This means that pure collocation will not allow for scaling to large problem sizes in terms of the number of patches used. If a least squares approach is instead used, the stability norm is not at all affected by the patch size. This is a very important property as it provides robustness for large numbers of patches. We can also note that the error curves for the collocation case are quite irregular. There are two main reasons. First, using a global node set means that the alignment of the nodes and the patches varies with the particular choices of $h$ and $H$, as illustrated by the oscillations in the error curves. Second, the global node set is Cartesian, which is sub-optimal for the RBF-QR method, because nodes on a grid are typically not unisolvent for polynomials
and this results in some numerical issues~\cite{LaFo05,FoLaFly11,LLHF13}.  

\begin{figure}[!htb]
\centering
\includegraphics[height=0.38\textwidth]{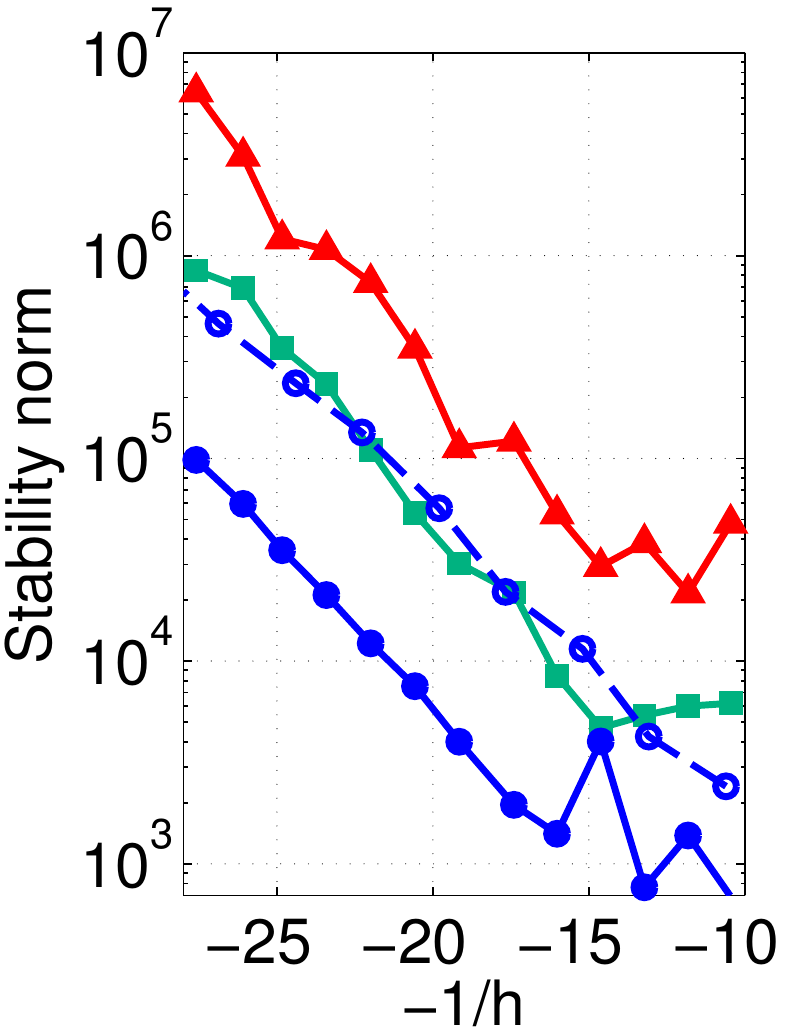}
\includegraphics[height=0.38\textwidth]{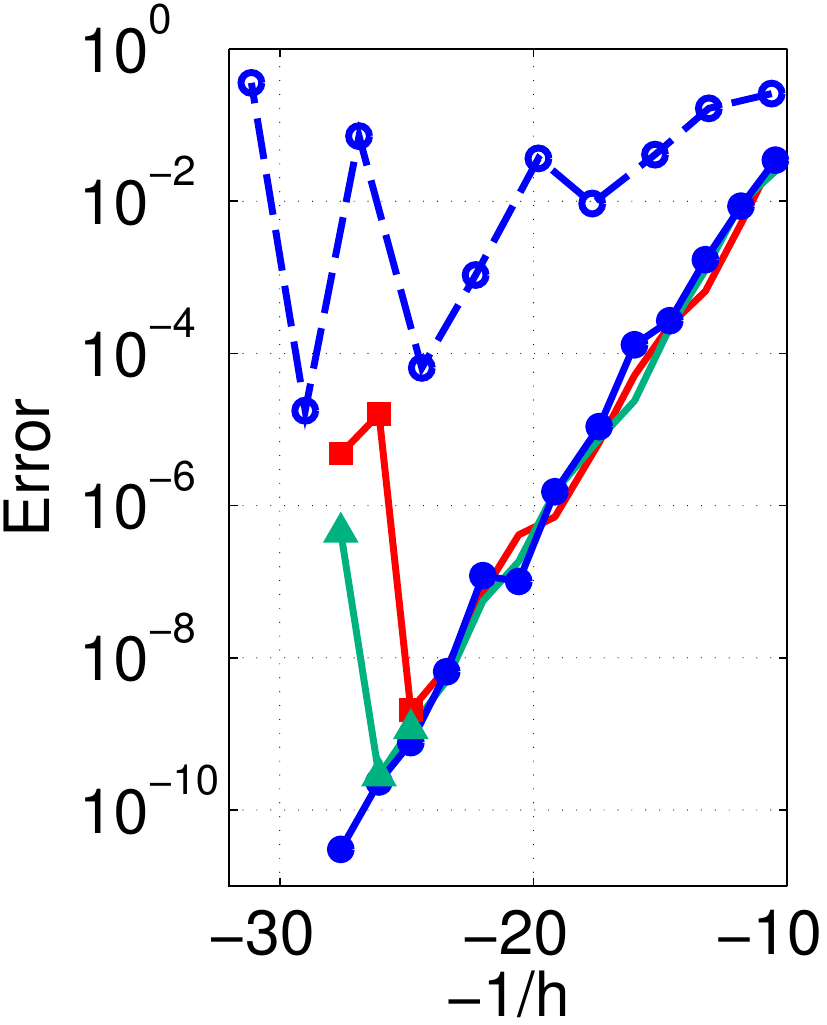}
\caption{The numerically estimated stability norm (left) and the corresponding error (right) as a function of $-1/h$ for oversampling $\beta=1.1$ ($\triangle$), $\beta=1.2$ ($\Box$), and $\beta=1.5$ ($\bigcirc$) for LS-RBF-PUM (solid lines, solid markers) and C-RBF-PUM (dashed lines, open markers).}
\label{fig:fillnorm}
\end{figure}

The second case we consider is fixed patch size $H$ and varying fill distance $h$. We test three different choices of oversampling for LS-RBF-PUM. The results are shown in Figure~\ref{fig:fillnorm}. The stability norm grows exponentially for both the collocation and the least squares approach. A higher degree of oversampling reduces the stability norm. It is worth to notice the effect this has on the error. The results for the three different cases of oversampling are very similar down to the points where the convergence trend is lost. The point of departure from the common trend can be approximately identified by multiplying the stability norm with the machine precision ($\sim 10^{-16}$) and comparing with the error, as predicted by the error estimate~\eqref{eq:errest}. The stability norm effectively captures the effective conditioning of the problem for LS-RBF-PUM. In the collocation case, the accuracy is not high enough for the stability norm to come into play. 

\begin{figure}[!htb]
\centering
\includegraphics[height=0.38\textwidth]{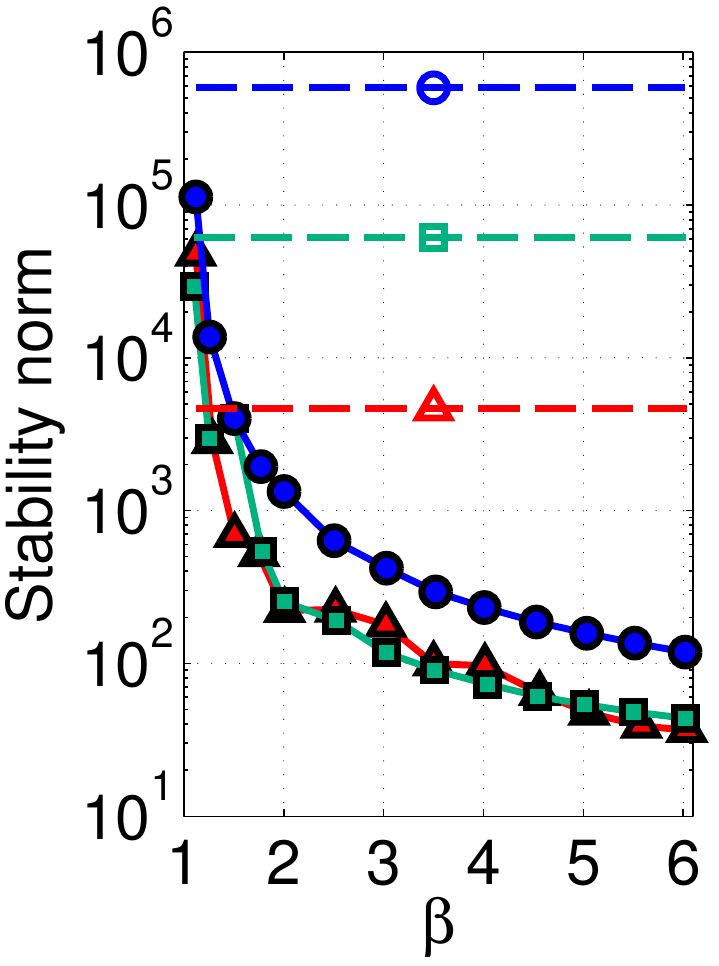}
\includegraphics[height=0.38\textwidth]{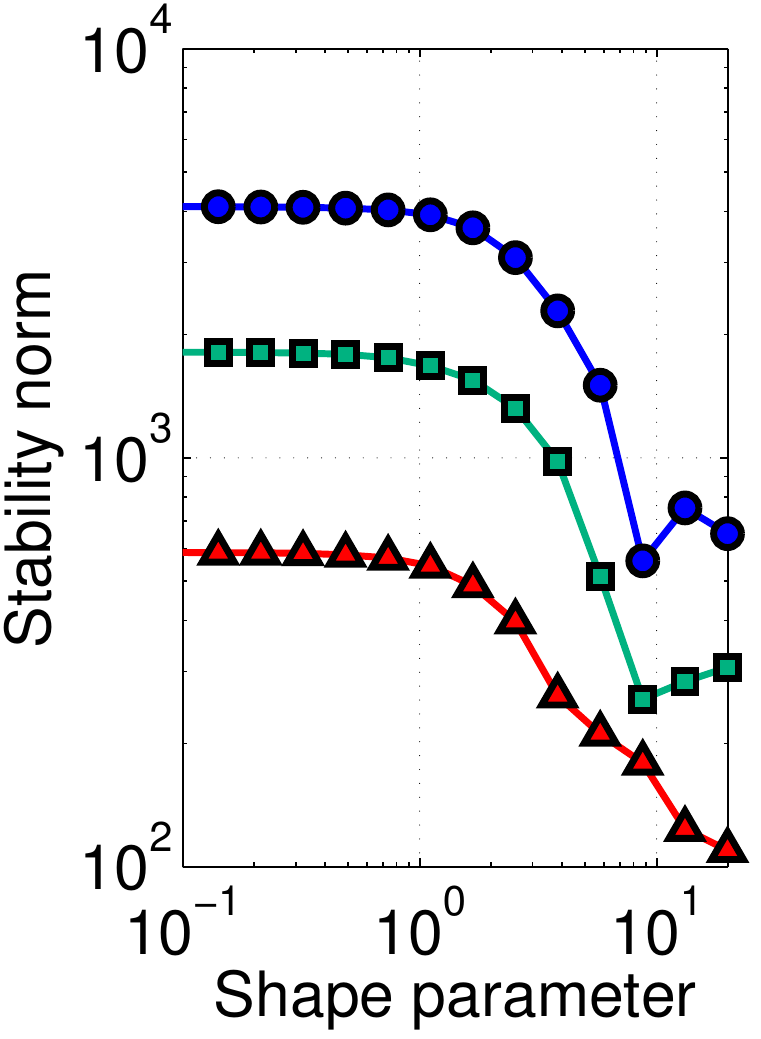}
\includegraphics[height=0.365\textwidth]{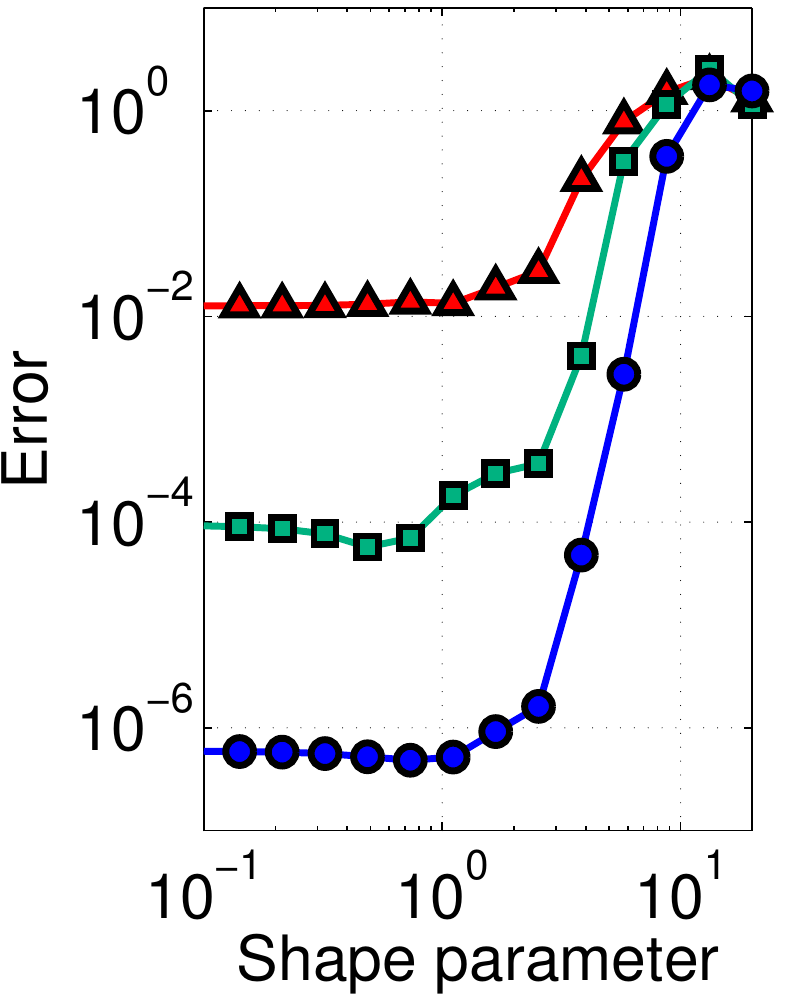}
\caption{The numerically estimated stability norm as a function of $\beta=M/N$ for $H=0.4$ ($P=100$) (left), the stability norm as a function of the shape parameter $\ep{}$ for $H=4/11$ ($P=121$) (middle), and the error as a function of $\ep{}$ for test function $u_2$ (right). Results are shown for LS-RBF-PUM (solid lines, solid markers), and in the left subfigure also for C-RBF-PUM (dashed lines, open markers), for $n=28$ ($\triangle$), $n=55$ ($\Box$), and $n=91$ ($\bigcirc$).}
\label{fig:shape}
\end{figure}

The amount of oversampling can be increased to reduce ill-conditioning, especially for larger numbers of points per patch. In the left part of Figure~\ref{fig:shape}, we investigate the relation between the stability norm and the rate of oversampling $\beta$ for fixed $H$ and $h$. The stability norm decreases rapidly initially and then levels out at a low level. The errors are mostly unaffected. This means that it is possible to have stability for high resolutions or large problem sizes by paying the computational price of having a larger oversampling rate.

Finally, the shape parameter, which has a crucial effect on conditioning when stable evaluation methods are not employed, is investigated here. RBF-QR provides stability for small values of $\ep{}$, and as can be seen in Figure~\ref{fig:shape}, the stability norm is constant in this regime. The stability norm decreases for increasing $\ep{}$, but for smooth functions, this is a regime where also the accuracy is lower.

\subsection{Computational efficiency}
\label{sec:eff}

\begin{figure}[!htb]
\centering
\includegraphics[height=0.38\textwidth]{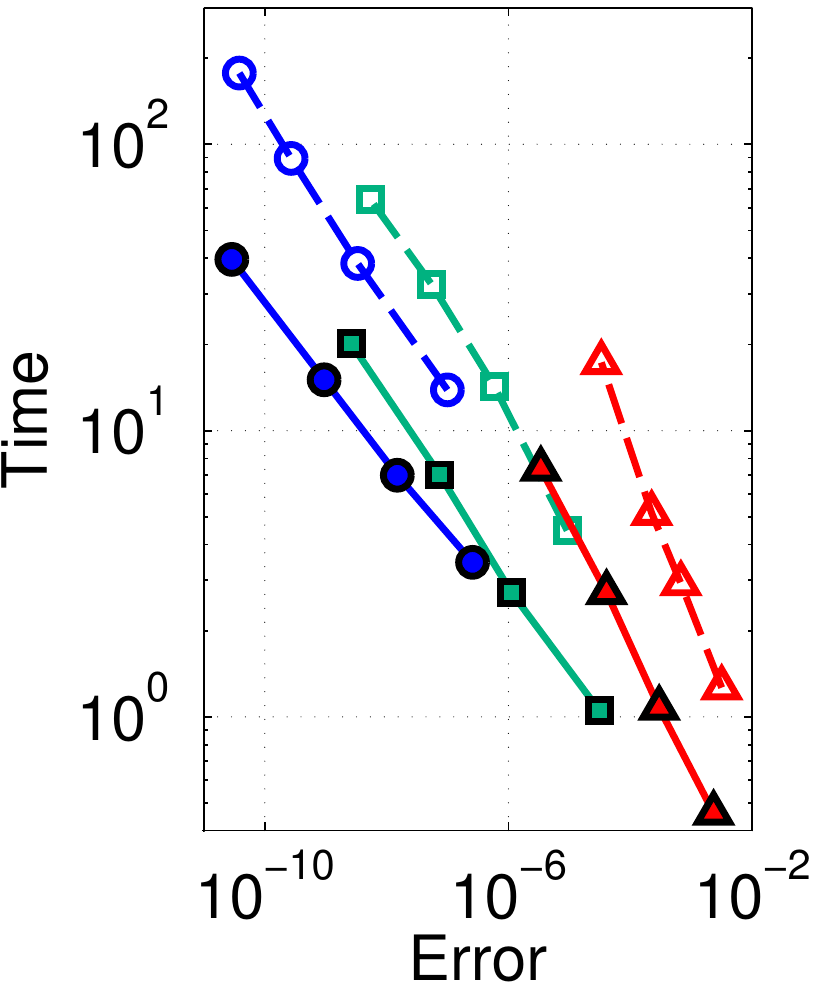}
\includegraphics[height=0.38\textwidth]{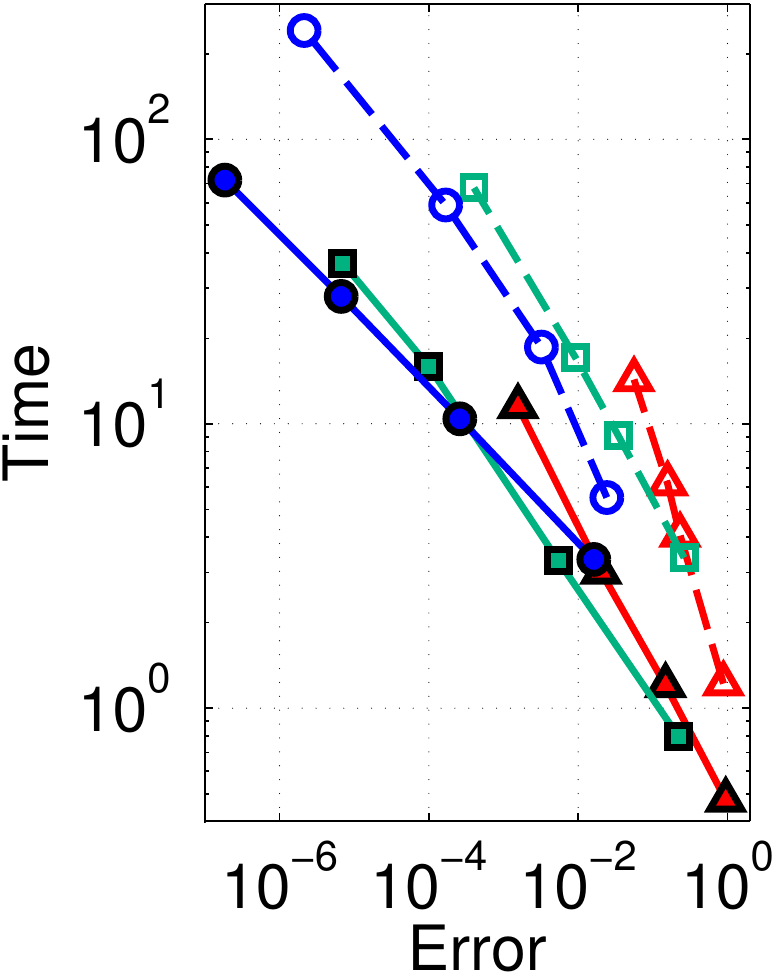}
\caption{Computational time against error for varying $H$ and $n=28$ ($\triangle$), $n=55$ ($\Box$), and $n=91$ ($\bigcirc$) for LS-RBF-PUM (solid lines, solid markers) and C-RBF-PUM (dashed lines, open markers) for the hyperbolic sine function $u_1$ (left) and the truncated sum $u_4$ (right). }
\label{fig:time}
\end{figure}

In Figure~\ref{fig:time}, we compare the computational efficiency of LS-RBF-PUM with that of C-RBF-PUM. 
LS-RBF-PUM is 5--10 times faster in all cases, and the gain increases with problem size. This is a combined effect of the reduced setup cost due to the identical local node layouts, the more efficient use of the degrees of freedom, and the increased robustness for larger problems.

With the possibility to vary both $H$ and $h$ it is possible to reach the same accuracy in different ways. The question is then which way is the most computationally efficient. There is no unique answer as can be seen in the figure. Rather, a smaller number of points per patch should be used if a low accuracy is desired, while more points per patch should be applied to reach a higher accuracy.

The level of accuracy that can be reached within a certain time depends on the function that is approximated. Clearly the function $u_4$ requires a higher resolution than $u_1$ for a given target accuracy, but the relation between the two methods and the refinement strategies are similar in both cases.

\subsection{Irregularly shaped domains}
\begin{figure}[!htb]
\centering
\includegraphics[height=0.415\textwidth]{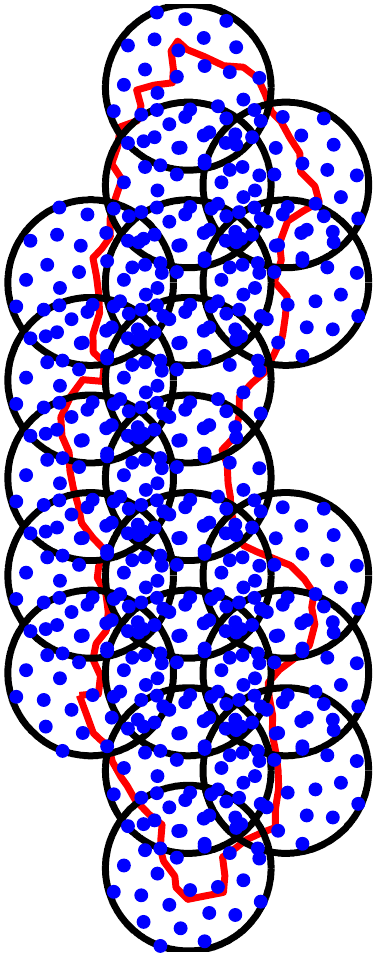}\hspace*{1cm}
\raisebox{2.5mm}{\includegraphics[height=0.38\textwidth]{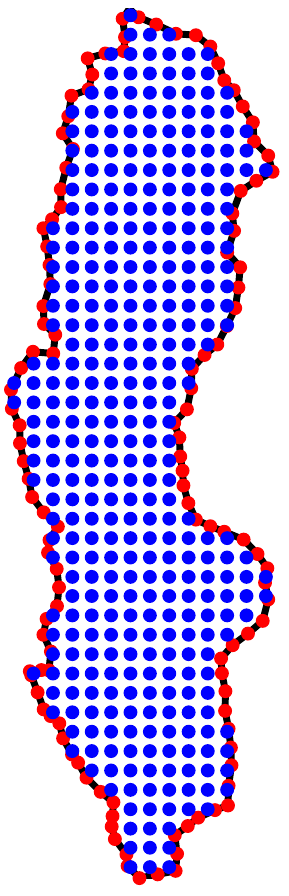}}
\caption{Patches with identically distributed local node sets covering the domain $\Omega_L$ (left), and least squares evaluation points on a Cartesian grid in the interior, and uniform with respect to arc length on the boundary (right) for LS-RBF-PUM.}
\label{fig:swe}
\end{figure}

As explained in the beginning of the section, the convergence experiments were performed on a square in order to promote regularity. Here, we verify that the results hold also on the irregularly shaped domains $\Omega_S$, see Figure~\ref{fig:nodes}, and $\Omega_L$, see Figure~\ref{fig:swe}.

\begin{figure}[!htb]
\centering
\includegraphics[height=0.365\textwidth]{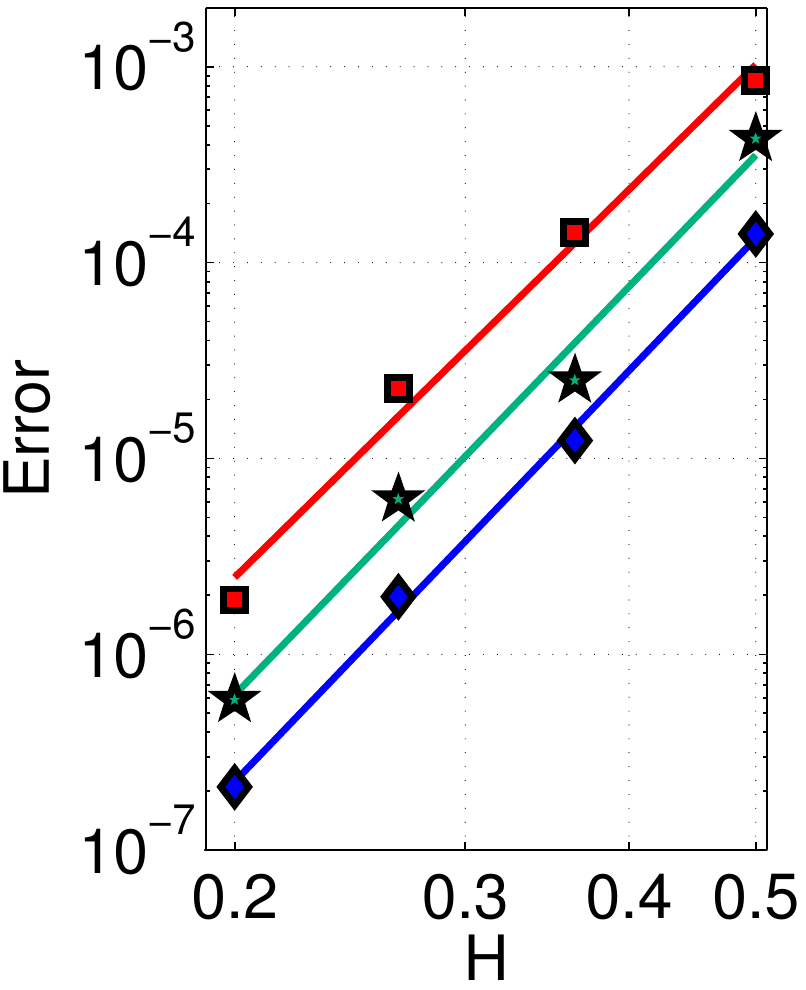}
\includegraphics[height=0.38\textwidth]{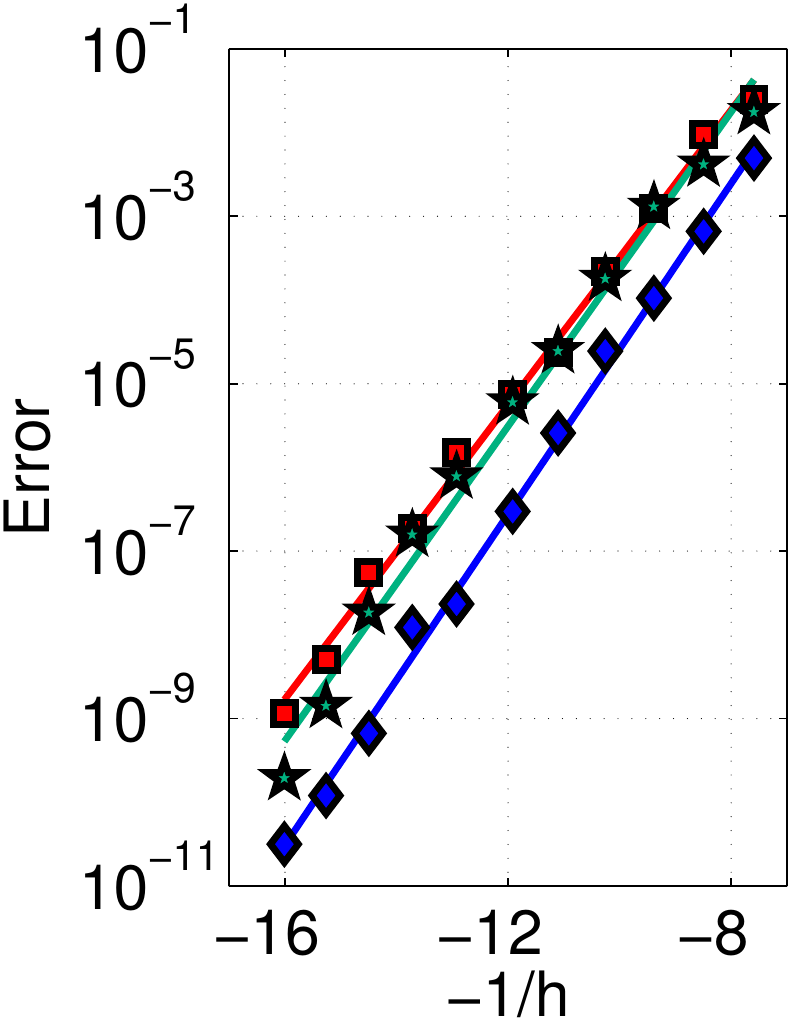}
\caption{Algebraic convergence results with $n=55$ (left) and exponential convergence results with $H=0.4$ (right) for the domains $\Omega_B$ ($\Box$), $\Omega_S$ ($\star$), and $\Omega_L$ ($\diamond$) using LS-RBF-PUM for the trigonometric solution function $u_2$.}
\label{fig:irreg}
\end{figure}

Results for the three domains for the algebraic and exponential convergence modes are shown in Figure~\ref{fig:irreg}. The convergence is regular for all three domains and behaves according to theory. The errors are smaller for the irregular domains. The explanation is that the square is larger (enclosing both of the other domains) and therefore contains more of the function that is approximated.

\subsection{Numerical convergence results for three-dimensional problems}

\begin{figure}[!htb]
\centering
\begin{tabular}{@{}cc@{}}
\includegraphics[height=0.36\textwidth]{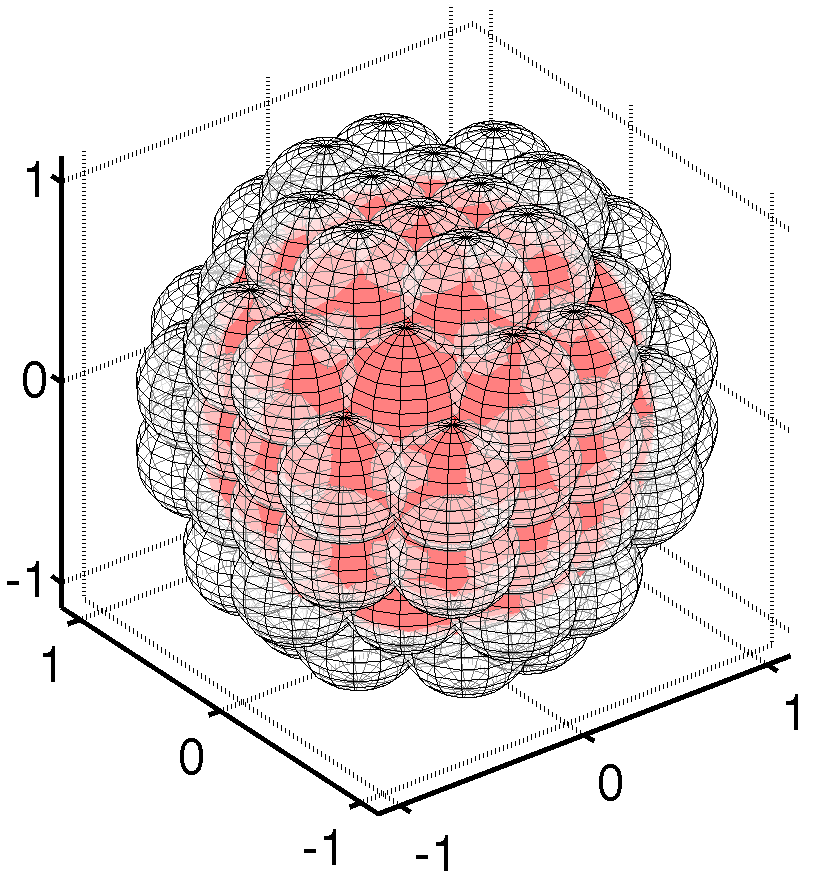} &
\includegraphics[height=0.42\textwidth]{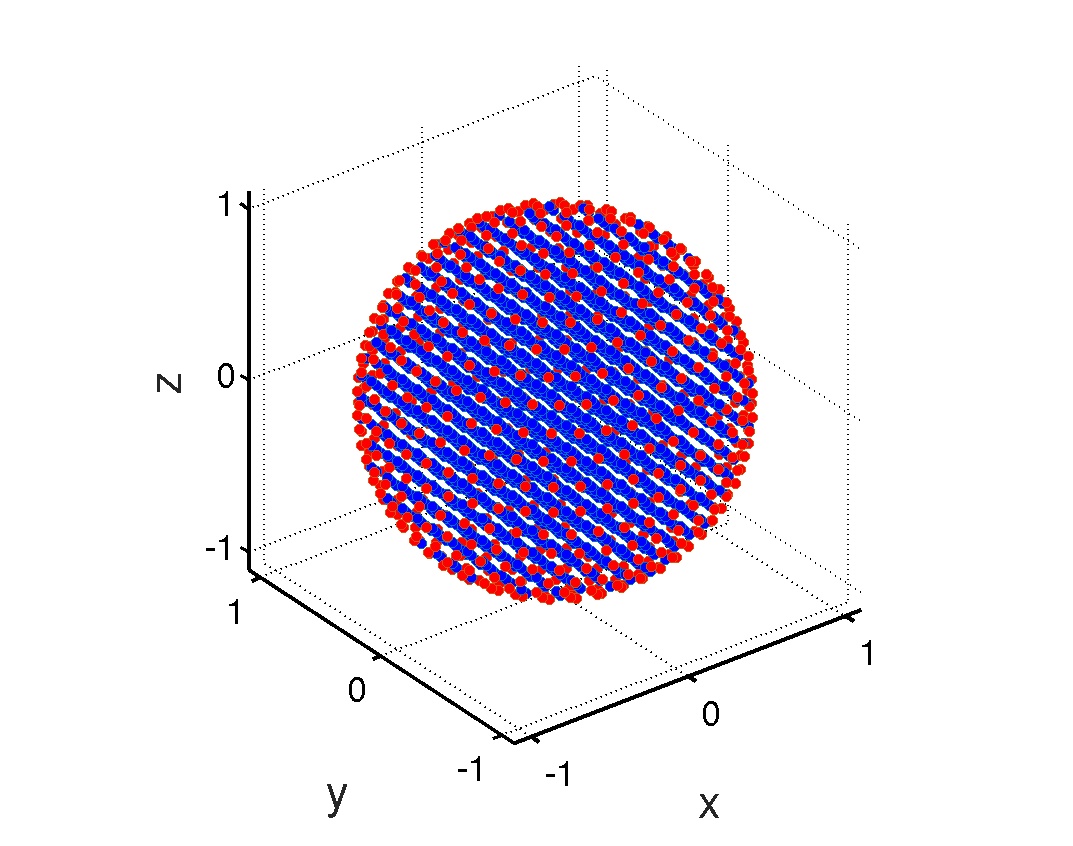}\\
\raisebox{6mm}{\hspace{1mm}\includegraphics[width=0.26\textwidth]{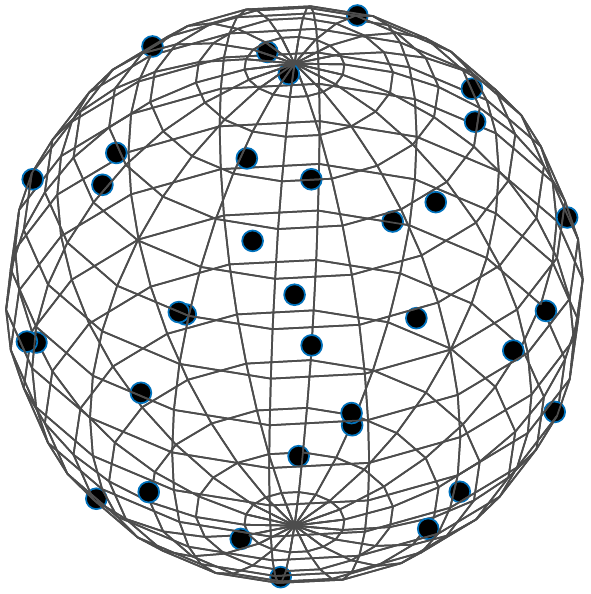}}&
\includegraphics[height=0.38\textwidth]{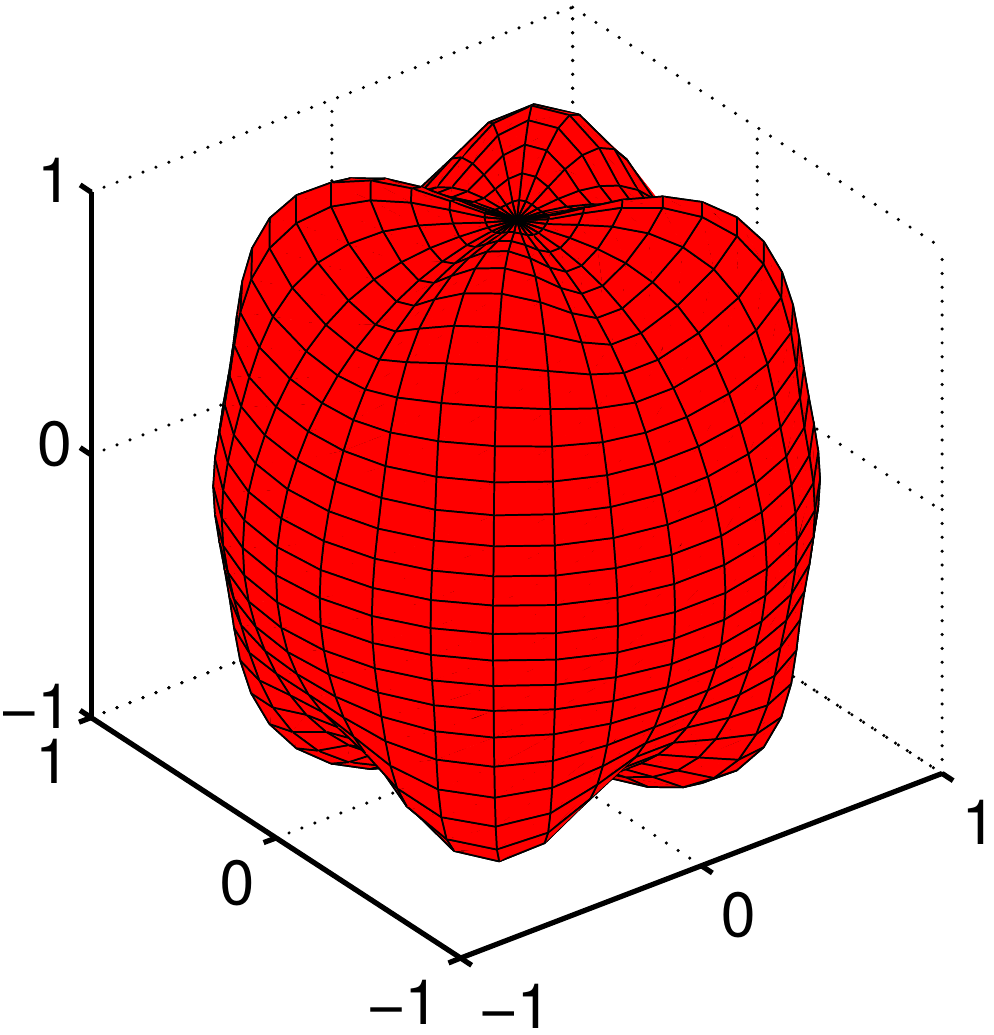} 
\end{tabular}
\caption{Patches on the spherical domain $\Omega_U$ (top left), quasi uniformly distributed least squares points inside $\Omega_U$ and on the surface $\partial\Omega_U$ (top right), a single patch (enlarged) with $n=35$ local node points (bottom left), and the star shaped domain $\Omega_Q$ (bottom right).}
\label{fig:dom3}
\end{figure}

For the three-dimensional test cases, we perform convergence experiments only for LS-RBF-PUM to confirm that the results are similar to the observations for the two-dimensional case. The geometry of the two domains used, $\Omega_U$ and $\Omega_Q$, with spherical patches, together with examples of local node sets and least squares evaluation points are displayed in Figure~\ref{fig:dom3}.

\begin{figure}[!htb]
\centering
\includegraphics[height=0.38\textwidth]{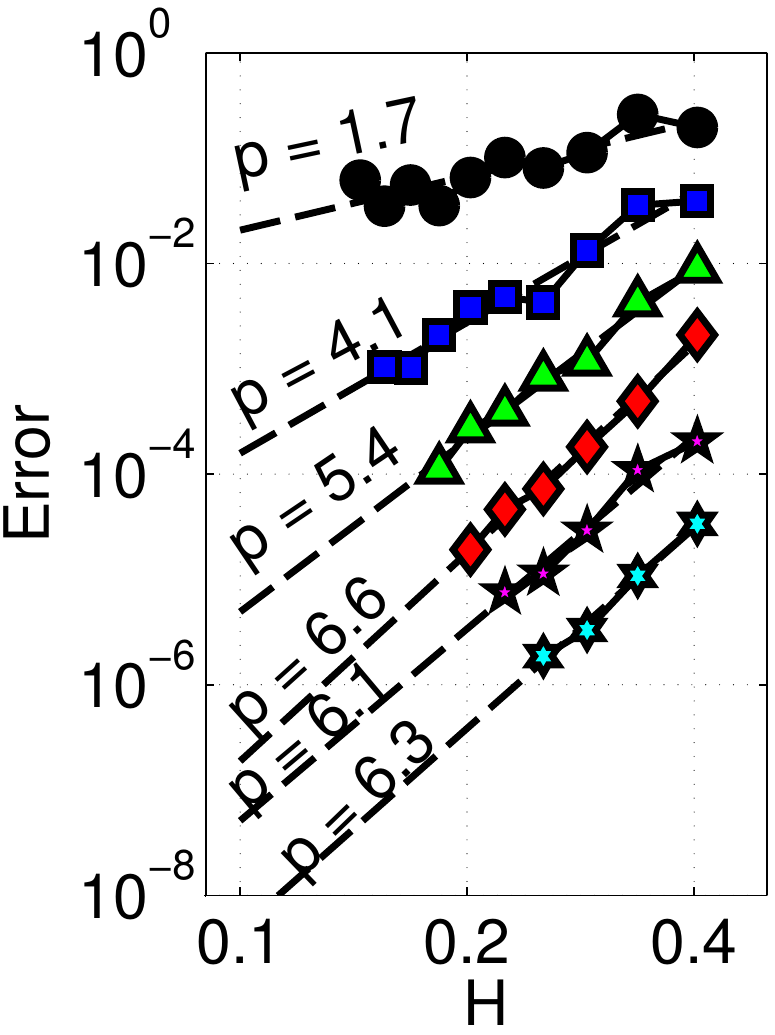}
\includegraphics[height=0.38\textwidth]{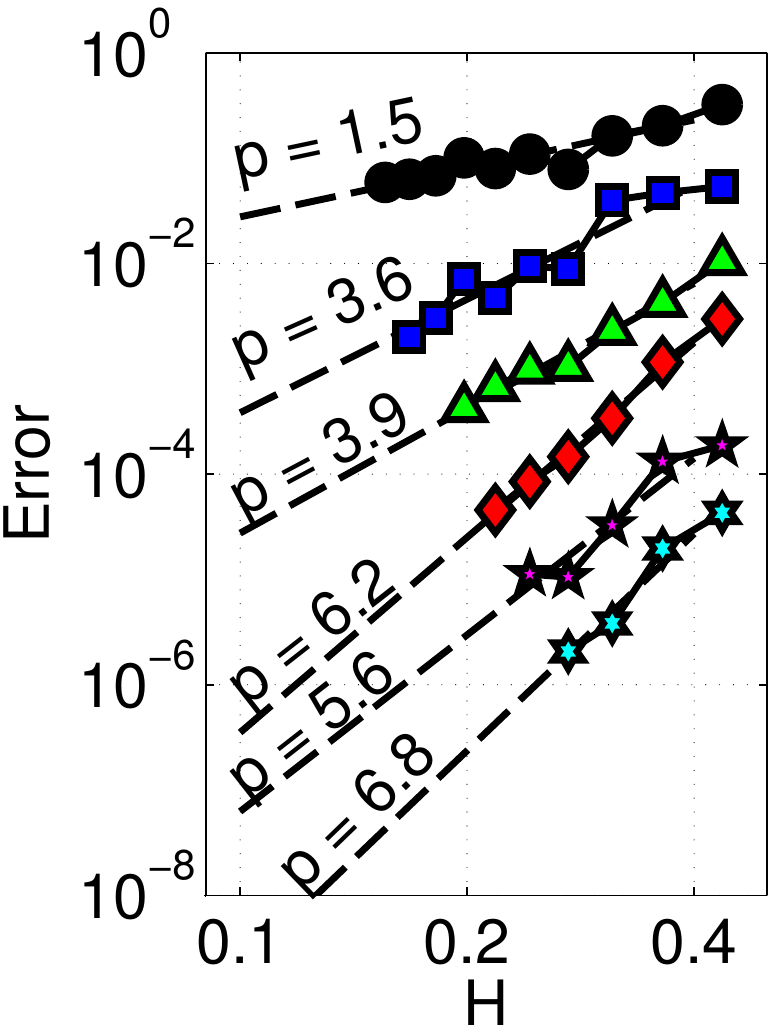}
\includegraphics[height=0.38\textwidth]{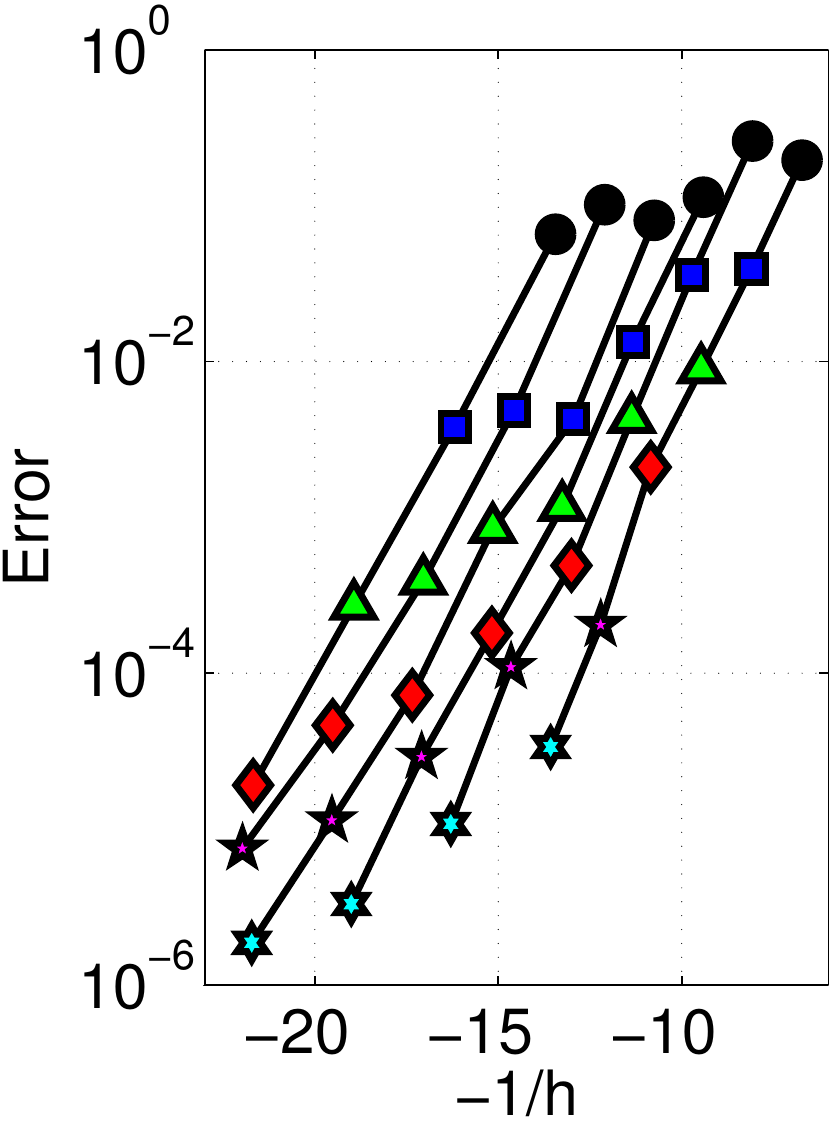}
\caption{Algebraic convergence results for the spherical domain $\Omega_U$ (left) and the star shaped domain $\Omega_Q$ (middle), and spectral convergence results for $\Omega_U$ (right). All results are for LS-RBF-PUM with $n=20$ ($\bigcirc$), $n=35$ ($\Box$), $n=56$ ($\triangle$), $n=84$ ($\diamond$), $n=120$ ($\star$), $n=165$ (hexagram). In the right subfigure, $H=2.02/k$ for $k=5,\ldots,10$, from right to left corresponding to $P=81$, $136$, $179$, $280$, $365$ and $551$. For the horizontal axis, $h$ represents the average fill distance.}
\label{fig:3dconv}
\end{figure}

The convergence results for fixed numbers of points per patch $n$ and varying patch size $H$ are shown in the left and middle subfigures of Figure~\ref{fig:3dconv}. Using Estimate~\ref{est:alg}, we would expect algebraic convergence rates of orders $0.5$, 1.5, 2.5, 3.5, 4.5, and 5.5 for $n=20$, 35, 56, 84, 120, 165. However, the numerically estimated rates are between 0.8 and 3.1 orders higher. The average difference between the numerically estimated order and the expected order is 1.8. Going back to~\cite{RieZwi10}, which provides the underlying estimates for RBF interpolation, we can see that in exchange for a larger constant, we can replace the term $-d/2$ in the convergence order with $-d/s$, where $1\leq s<\infty$. That is, the estimate allows for an improvement of up to almost $d/2=1.5$ in the order of convergence, depending on which norms are used in the underlying estimate. 

The right subfigure of Figure~\ref{fig:3dconv} shows convergence as a function of fill distance $h$ for different values of $H$. In the plot, we have not used the precise fill distance according to~\eqref{eq:filldist} due to irregularities in the node sets. Instead we have approximated the average fill distance through $h\approx H/n^{1/3}$. It is clear from the figure that for a given $h$, the accuracy is improved as $H$ is increased. This can be understood from the fact that a global approximation is the most accurate way to represent a function given a certain node density. However, as discussed in Section~\ref{sec:eff} it is not the best choice from the computational efficiency point of view. Furthermore, as shown in Section~\ref{sec:stab}, the ill-conditioning of a problem increases rapidly with $n\approx(H/h)^d$. The slopes are similar for all values of $H$ as predicted by Estimate~\ref{est:spec}, where the rate $\gamma$ is independent of $H$. 
%
%
Qualitatively, the convergence results in two and three dimensions agree with each other as well as with the theoretical results.

\section{Discussion}
In this paper, we have proposed a new least squares formulation of RBF-PUM. Just as  C-RBF-PUM (the original collocation based formulation), LS-RBF-PUM requires the RBF-QR method or another stable evaluation method
in order to converge as the patch size is refined. RBF-QR~\cite{FoPi07,FoLaFly11,LLHF13} is currently only available in up to three space dimensions. However, from experience an accuracy of about $10^{-5}$ can be achieved without a stable method, which is often enough for practical purposes, especially when working with high-dimensional problems.

LS-RBF-PUM significantly simplifies the handling of geometry. The node points and patches do not need to conform to the geometry, and the method is not sensitive to the location of the least squares evaluation points relative to the geometry. High quality node points can even be pre-computed and stored, since the patch geometry only depends on the dimension.



In this paper, we have derived the first theoretical estimates for RBF-PUM solutions to PDEs. The numerical results show that the actual error behavior can be understood from the theoretical results. The matrix norm that appears in the estimates has been investigated numerically. The most important results from a practical perspective are (i) that the norm does not grow at all under patch refinement, which means that LS-RBF-PUM can be used for solving large scale problems, and (ii) that by increasing the amount of oversampling $\beta$, the norm can be made small, which means that we can reduce the conditioning to reach a higher accuracy if we are willing to pay the added computational cost.

Even though LS-RBF-PUM uses more node points than C-RBF-PUM for the same spatial resolution, LS-RBF-PUM is 5-10 times faster than C-RBF-PUM. The main reasons for the gain are the decreased setup cost because local node points are identically distributed with respect to the patches, and the more efficient use of the degrees of freedom when $n$ is the same in all patches. For C-RBF-PUM, $n$ becomes smaller in boundary patches. This can be overcome by making boundary patches larger~\cite{LaHer17}, but then makes the algorithm more complicated.


An improvement that has not been investigated here is to make the approximation adaptive. This is done for interpolation with good results with respect to accuracy in~\cite{CaDeRoPe16b}, and would be highly relevant in the PDE context.


\bibliographystyle{siam}
\bibliography{refs}

\end{document}